\documentclass[12pt]{elsarticle}
\usepackage{amsmath,amssymb,epsfig,amsfonts,epsfig,graphicx}
\usepackage{stmaryrd} 
\SetSymbolFont{stmry}{bold}{U}{stmry}{m}{n} 
\usepackage{epsfig}
\usepackage{tikz-cd}
\usepackage{color}
\usepackage{pict2e}
\usepackage{bbold}

\journal{SIAM Journal on Mathematical Analysis}

\newtheorem{Remark}{Remark}
\newtheorem{Lemma}{Lemma}
\newtheorem{Definition}{Definition}

\newtheorem{Theorem}{Theorem}
\newtheorem{Proposition}{Proposition}

\renewcommand{\rho}{\varrho}
\renewcommand{\d}{\mathrm{d}}
\newcommand{\ds}{\,\mathrm{d}s}

\newcommand{\dx}{\,\mathrm{d}x}

\newcommand{\emin}{\underline{e}_{\mu,\mu_c}}
\newcommand{\wmin}{\underline{w}_{\mu,\mu_c}}
\newcommand{\muccrit}{\mu_c^{\rm crit}}
\newcommand{\Zmu}{Z_{\mu,\mu_c}}
\newcommand{\Smu}{S_{\mu,\mu_c}(\alpha)}
\newcommand{\dy}{\,\mathrm{d}y}

\newcommand{\ve}{\varepsilon}
\newcommand{\sa}{\sin(\alpha)}
\newcommand{\ca}{\cos(\alpha)}
\newcommand{\ao}{\widetilde{\alpha}_\mathrm{opt}}

\newcommand{\aen}{\alpha_{\varepsilon_n}}
\newcommand{\uen}{u_{\varepsilon_n}}
\newcommand{\ta}{\widetilde{\alpha}}
\newcommand{\tu}{\widetilde{u}_{\ve_n}}
\newcommand{\aopt}{\alpha_\mathrm{opt}}

\newcommand{\shalf}{\sin^2\!\big(\frac{\alpha}{2}\big)}
\newcommand{\Wred}{\widehat{W}}
\newcommand{\Ered}{\widehat{E}}
\renewcommand{\AA}{\alpha_1^-}
\newcommand{\BB}{\alpha_1^+}

\newcommand{\cored}{\widehat{c}_0}
\newcommand{\uhom}{\overline{u}}
\newcommand{\nn}{\nonumber}

\newcommand{\R}{\mathbb{R}}
\newcommand{\N}{\mathbb N}

\newcommand{\tr}{\mathrm{tr}}
\newcommand{\SO}{\mathrm{SO}}

\newcommand{\DR}{\text{D}}
\newcommand{\Curl}{\text{Curl}}
\newcommand{\dev}{\text{dev}}
\newcommand{\sym}{\text{sym}}
\newcommand{\id}{\mathbb{1}}

\begin{document}

\begin{frontmatter}
\title{The Gamma-limit of the simple shear problem\\
in nonlinear Cosserat elasticity}
\author{Thomas Blesgen}
\ead[1]{t.blesgen@th-bingen.de}
\address{Bingen University of Applied Sciences, Berlinstra{\ss}e 109,
D-55411 Bingen, Germany}
\author{Patrizio Neff}
\ead[2]{patrizio.neff@uni-due.de}
\address{University of Duisburg-Essen, Faculty of Mathematics,
Thea-Leymann-Stra{\ss}e 9, D-45127 Essen, Germany}
\cortext[2]{Corresponding author}
\cormark[2]

\begin{abstract}
The zero and first order Gamma-limit of vanishing internal length scale are
studied for the mechanical energy of a shear problem in geometrically nonlinear
Cosserat elasticity. The convergence of the minimizers is shown and the limit
functionals are characterized.
\end{abstract}

\end{frontmatter}

\section{The Cosserat model in simple shear}
\label{secintro}
\begin{figure}[h!tb]
\unitlength1cm
\begin{picture}(11.5,4.0)
\put(4.5,0.0){\psfig{figure=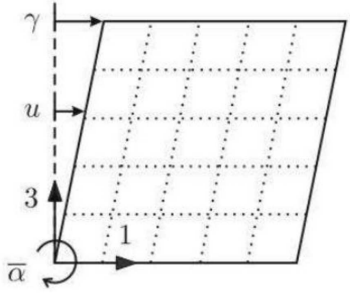,width=6.3cm}}
\end{picture}
\caption{\label{fig1}
\footnotesize The deformed state exhibits a homogeneous region in the
interior of the structure which motivates the kinematics of simple shear.}
\end{figure}

We investigate the deformation of an infinite layer of material in 3D with unit
height, fixed at the bottom and sheared in $e_1$-direction with amount
$0<\gamma<2$ at the upper face, cf. Fig~\ref{fig1}.
Within a geometrically non-linear Cosserat theory,
\cite{Capriz89,Cos1909,Cos1991}, the mechanical behaviour of
the material can be modelled with the help of the standard deformation map
$\varphi:\widehat{\Omega}\to\R^3$ and the tensor field of orthogonal
micro-rotations $R:\widehat{\Omega}\to\SO(3)$, describing the translation and
independent rotation of a material point, respectively. Here,
$\widehat{\Omega}\!\subset\!\R^3$ is the reference configuration.

By $\mu>0$ we denote the standard elastic shear modulus, $\mu_c\ge0$ is the
Cosserat couple modulus, $\lambda\in\R$ the second elastic Lam{\'e} parameter,
$a_1\ge0$, $a_2\ge0$, $a_3\ge0$ are non-dimensional constants; $L_c>0$
is the characteristic length scale.

The deformed material is characterized by the minimizers $(\varphi,R)$ of the
isotropic mechanical energy
{\small
\begin{align}
E_{3D}(\varphi,R) &= \int_{\widehat{\Omega}}\mu\big|
\sym(\overline{U}-\id_3)\big|^2
+\mu_c\big|\text{skew}(\overline{U}-\id_3)\big|^2
+\frac{\lambda}{4}\Big[(\det\overline{U}-1)^2+\Big(
\frac{1}{\det\overline{U}}-1\Big)^2\Big]\nn\\
&\quad\quad +\mu\frac{L_c^2}{2}\Big(a_1\big|\dev\, \sym\, R^T\Curl R\,
\big|^2+a_2\big|\text{skew}R^T\Curl R\big|^2+\frac{a_3}{3}\tr(R^T\Curl R)^2
\Big)\dx\nn\\
\label{E3d}
&= \int_{\widehat{\Omega}}W_{\text{mp}}(\overline{U})+W_{\text{disloc}}
(R^T\Curl R)\dx\,,\qquad \overline{U}=R^T\DR\varphi
\end{align}
}
\noindent subject to certain boundary conditions, see
\cite{BN22,Neff15,Neff06,Neff08,Neff09,Ble14,Ble13} for
further information and \cite{Dmit09} for a comparison to experiments.

The symmetry of the boundary conditions and the infinite extension in
$e_1$-direction lead to the reduced kinematics
\begin{align}
\label{Fdef}
\varphi(x_1,x_2,x_3)=\begin{pmatrix}
x_1+u(x_3)\\x_2\\x_3
\end{pmatrix}\,,\quad
F=\DR\varphi(x_1,x_2,x_3)=\begin{pmatrix} 1&0&u'(x_3)\\0&1&0\\0&0&1
\end{pmatrix}\,,
\end{align}
with $u(0)=0$ and $u(1)=\gamma$.
The microrotations $R\in\SO(3)$ satisfy the identity
\begin{align}
\label{Rdef}
R(x_1,x_2,x_3)=\begin{pmatrix}
\cos\alpha(x_3)&0&\sin\alpha(x_3)\\0&1&0\\
-\sin\alpha(x_3)&0&\cos\alpha(x_3) \end{pmatrix}
\end{align}
with fixed axis of rotation $e_2$, implying
\begin{equation}
\label{curlid}
\Curl R=\begin{pmatrix}
0 & -\sin\alpha(x_3)\alpha'(x_3) & 0\\ 0 & 0 &0\\
0 & -\cos\alpha(x_3)\alpha'(x_3) & 0 \end{pmatrix}\,.
\end{equation}
From now on, we denote $x_3$ by $x$ and set
\[ \Omega:=(0,1)\subset\R. \]

The deformation $u$ introduced in (\ref{Fdef}) and the local rotation angle
$\alpha$ around the $x=x_3$-axis in (\ref{Rdef}) permit to re-write the
functional $E_{3D}$ in a simpler form, see \cite{BN22}.
The solutions to the 3D-shear problem can thus be obtained as minimizers of the
{\it mechanical energy functional} and the
{\it reduced mechanical energy functional}
\begin{align}
E(u,\alpha) \;=\;\; &\frac{\mu}{2}\int_0^14L_c^2|\alpha'|^2+|u'|^2
+\Big(\sa u'-4\shalf\!\Big)^2\dx\nn\\
\label{Edef}
& +\frac{\mu_c}{2}\int_0^1\Big(\ca u'-2\sa\Big)^2\dx,\\
\label{Ered}
\Ered(u,\alpha) \;=\; &\frac{\mu}{2}\int_0^1\!4L_c^2|\alpha'|^2+|u'|^2
+[\alpha(\alpha\!-\!u')]^2\dx+\frac{\mu_c}{2}\!\int_0^1\!\Big(
\frac{2\!-\!\alpha^2}{2}u'-\frac{6\alpha\!-\!\alpha^3}{3}\Big)^2\dx.
\end{align}

The functional $\Ered(u,\alpha)$ is obtained from $E(u,\alpha)$ in
(\ref{Edef}) after introducing the third-order expansions
$\ca\sim1-\frac{\alpha^2}{2}$,
$\sa\sim\alpha-\frac{\alpha^3}{6}$ and dropping all higher order terms
except $\frac{\mu_c}{2}\Big(\frac{\alpha^4}{4}|u'|^2-\frac13\alpha^5u'
+\frac19\alpha^6\Big)$ to finally get a quadratic form.

\vspace*{2mm}
In this article we are concerned with the zero-order and first-order
Gamma-limit of vanishing internal length scale $L_c$ of $E(u,\alpha)$.
However, the methods of this article are also applicable to the Gamma-limit
$L_c\searrow0$ of $\Ered(u,\alpha)$ and similar results can be obtained.

\vspace*{2mm}
Formally writing $L_c=\frac{\ve}{\sqrt{2\mu}}$ for $\ve\ge0$ leads to
\begin{align}
\label{Eeps}
E_\ve(u,\alpha) & \;:=\; \left\{ \begin{array}{ll}
\int_0^1\ve^2|\alpha'|^2+\frac{\mu}{2}|u'|^2
+W(u',\alpha)\dx, & \mbox{ if }(u,\alpha)\in(W^{1,2}(\Omega))^2,\\ +\infty, &
\mbox{ else} \end{array}\right.
\end{align}
with the potential
\begin{align}
\label{Wdef}
W(u',\alpha) &\;:=\; \frac{\mu}{2}
\Big(\sa u'-4\shalf\!\Big)^2+\frac{\mu_c}{2}\Big(\ca u'-2\sa\Big)^2.
\end{align}
By $W^{m,2}(\Omega)$ we denote the Sobolev space of $m$-times weakly
differentiable functions in $L^2(\Omega)$.
Let $C_\#^\infty(\Omega;\R)$ denote the smooth functions $g:\Omega\to\R$ with
$g(0)=g(1)$ and let $W^{1,2}_\#(\Omega;\,\R)$ be the closure of
$C_\#^\infty(\Omega;\,\R)$ with respect to the $W^{1,2}$-norm, i.e. the
Sobolev functions $g\in W^{1,2}(\Omega)$ with identical traces at the boundary.
The minimization of $E$ in (\ref{Edef}) is carried out in the reflexive Banach
space
\begin{align}
\label{Xdef}
{\cal X} &:= {\cal X}_u\times{\cal X}_{\alpha}:=\Big\{u\!\in\!W^{1,2}(\Omega;
\,\R)\,\Big|\,u(0)\!=\!0,\,u(1)\!=\!\gamma\Big\}\!\times\!\Big\{
\alpha\!\in\!W_\#^{1,2}(\Omega;\,[0,2\pi])\Big\}.
\end{align}

The concept of $\Gamma$-convergence describes the asymptotic
behaviour of a family of minimization problems. It is arguably the most
natural way to study the convergence of variational problems as it supplies
information not only of the minimizers itself, but also of the convergence of
the variational problems. Theorem~\ref{theo1} below states the circumstances.

The characterization of the zero-order and first-order $\Gamma$-limit of
$E_\ve$ allows to qualitatively and quantitatively understand the model for
small characteristic length scale $L_c$ which would otherwise demand
simulations with ultra-high spatial resolution.
Clearly, the $\Gamma$-limit differs from the direct limit. This is illustrated
in the following non-commutative diagram.

\[ \begin{tikzcd}
E_\ve \arrow{r}{\ve\to0}
\arrow[swap]{d}{\Gamma\!-\!\lim\limits_{\hspace*{-10pt}\ve\to0}}
& \widetilde{E}_0 \arrow[equal, "/" marking]{dl}{\not=}\\
E_0 & \quad \end{tikzcd} \]

The $\Gamma$-limit functional $E_0$ of $E_\ve$ will be identified in
Proposition~\ref{prop1}. For comparison, the limit functional $\widetilde{E}_0$
for $L_c\searrow0$ differs in general from $E_0$ and is simply
(cf. Eqn.~(\ref{Edef}))
\begin{equation}
\label{E0tdef}
\widetilde{E}_0(u,\alpha):=\frac{\mu}{2}\int_0^1|u'|^2+\Big(\!\sin(\alpha)u'
\!-\!4\sin^2\!\big(\frac{\alpha}{2}\big)\!\Big)^2\dx+\frac{\mu_c}{2}\int_0^1
\!\Big(\!\cos(\alpha)u'\!-\!2\sin(\alpha)\!\Big)^2\dx.
\vspace*{-4mm}
\end{equation}
The condensed energy
$E^\mathrm{cond}(u):=\min_\alpha\widetilde{E}_0(u,\alpha)$
can be determined explicitly, see the appendix.
Accordingly, the limit $L_c\searrow0$ of the full functional in 3D is
\begin{equation}
\label{E3ddef}
\widetilde{E}_{3D}(\varphi,R)=\int_{\widehat{\Omega}}\mu\big|
\sym(\overline{U}\!-\!\id_3)\big|^2
+\mu_c\big|\text{skew}(\overline{U}\!-\!\id_3)\big|^2
+\frac{\lambda}{4}\Big[(\det\overline{U}\!-\!1)^2+\Big(
\frac{1}{\det\overline{U}}\!-\!1\Big)^2\Big]\dx
\vspace*{-4mm}
\end{equation}
and the condensed energy
$E_{3D}^\mathrm{cond}(\varphi):=\min_{R\in\SO(3)}\widetilde{E}_{3D}(\varphi,R)$
can also be determined explictly.
At this point, we are unable to compute the $\Gamma$-limit $L_c\searrow0$ of
the full three-dimensional functional $\widetilde{E}_{3D}$ which is why
we restrict ourselves here to the analysis of the shear problem.

The paper is organized in the following way. In Section~\ref{sectheory} we
recall the theory of $\Gamma$-convergence as needed later.
Section~\ref{seczero} deals with the zero-order $\Gamma$-limit of $E$ as
$L_c\searrow0$. Here also the minimizers of $E$ are classified depending
on the values of $\mu$, $\mu_c$ and the amount of shear $\gamma$.
In Section~\ref{secfirst}, the first-order $\Gamma$-limit of $E$ is computed.
To that end, the energy is rescaled first.
In the appendix we compare the zero-order $\Gamma$-limit with a direct
minimization of $E$ for $L_c=0$.

\begin{figure}[h!t]
\unitlength1cm
\begin{picture}(11.5,4.7)
\put(0.5,-0.45){\psfig{figure=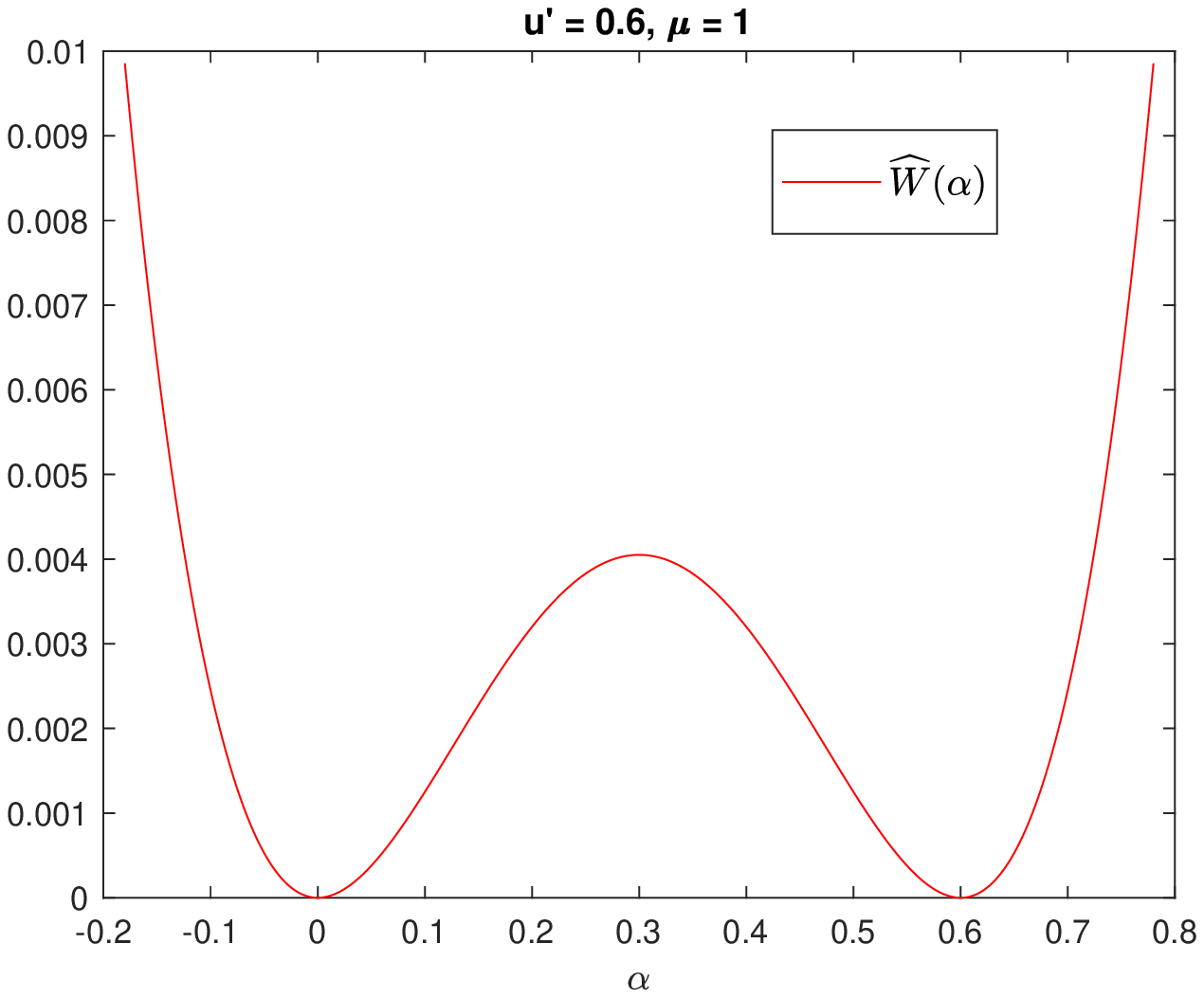,width=7.5cm}}
\put(8.0,-0.45){\psfig{figure=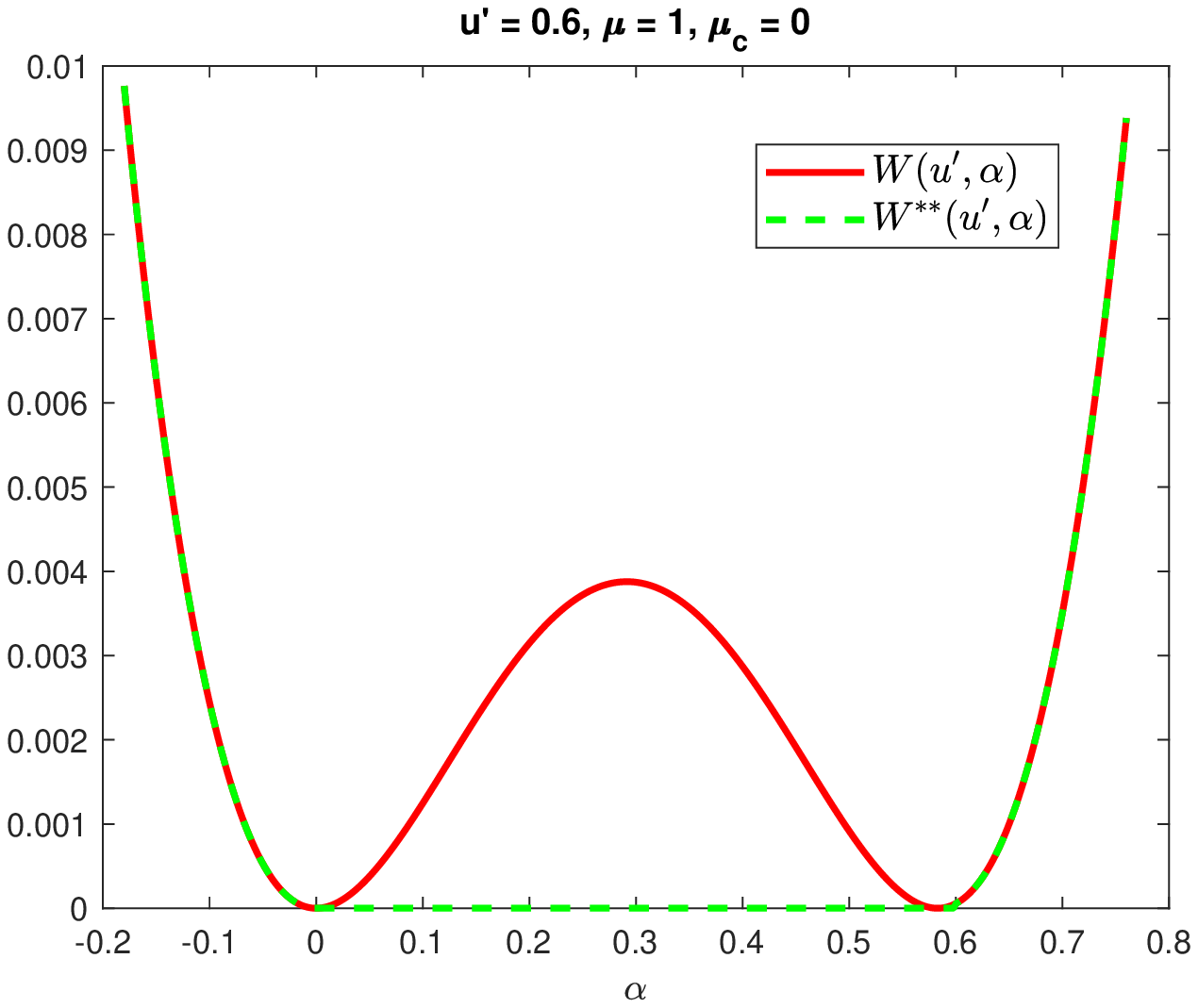,width=7.5cm}}
\end{picture}
\caption{\label{fig2}
Comparison of the reduced energy
$\Wred(\alpha)=\frac\mu2\big[\alpha(\alpha\!-\!u')\big]^2$ (left) and the full
energy $W$ (right) together with its convexification $W^{**}$ for
$\mu=1$, $\mu_c=0$, $u'=0.6$ and $\alpha\in[-0.2,0.8]$.}
\end{figure}

\vspace*{1mm}
As a good starting point and for gaining first insights into the concepts of
this article, we consider for fixed $u'$ the term
$\Wred(\alpha):=\frac{\mu}{2}\big[\alpha(\alpha\!-\!u')\big]^2$ from
Eqn.~(\ref{Ered}), see Fig.~\ref{fig2}. It is a double-well potential with
minima at $\alpha=0$ and $\alpha=u'$ and a maximum at $\alpha=\frac{u'}{2}$.
Similar potentials have been used for a long time to model phase transitions
and segmentation phenomena, see, e.g., \cite{CH86}.
(Recent articles on phase separation commonly replace the quartic polynomial by
a logarithmic expression closer to the correct physical free energy).

\vspace*{-2mm}
In the region separating the two minima, $\Wred$ is strictly concave, thereby
forming an energy barrier, also called surface energy. This is the minimal
amount of energy that must be provided to let the physical system pass between
the optimal states $\alpha=0$ and $\alpha=u'$. A related concept is
the activation energy of chemical reactions in an Arrhenius type equation,
see, e.g., \cite{Levine05}.

The energy barrier models analytically the resilience of the material to
changes of the inner molecular structure
needed to pass from the state $\alpha=0$ to the state $\alpha=u'$.
The associated molecular or atomistic restructuring
can mathematically be formulated as Markov processes, see, e.g., \cite{Lecca13}.

The mathematically correct form of the surface energy for the Cosserat
functional~(\ref{Edef}) will be analyzed in Section~\ref{secsurf}.

Since $\Wred$ originates from $W$ by a third-order Taylor expansion, one may
expect that also $W$ displays a double-well. Indeed, within a certain range of
$\mu$ and $\mu_c$, this is the case. Fig.~\ref{fig2} compares $W$ and
$\Wred$ for one set of parameters with striking similarity.

\vspace*{2mm}
If a double-well structure of $W$ and $\Wred$ is present, minimizing sequences
develop fine scale oscillations. For further explanations, we refer to
\cite{SM98} and references therein.

\vspace*{2mm}
The zero-order $\Gamma$-limit goes along with a convexification of $W$.
This convexification connects the minima and removes the energy barrier.
Using this convexified energy $E_0$, cf. Eqn.~(\ref{E0def}), should lead to a
significant improvement in the numerical simulations.

Fig.~\ref{fig2} shows only a part of $W$. Interestingly, the double-well is
very flat and located in a tiny section of the full graph of $W$ that may
easily be overlooked, cf. Fig.~\ref{fig4}.

\section{Theory of Gamma-convergence}
\label{sectheory}
\begin{Definition}
\label{def1}
Let $X$ be a topological space. A family of functionals
$(G_n)_{n\in\N}: X\to\overline{\R}:=\R\cup\{-\infty,+\infty\}$
{\em converges in the $\Gamma$-sense} for $n\to\infty$ to
$G: X\to\overline{\R}$, if the following two conditions are met\/:
\vspace*{-2mm}
\begin{itemize}
\item[(i)] (Liminf inequality)\\
For all $x\in X$ and every sequence $(x_n)_{n\in\N}\subset X$ with
$x_n\to x$ in $X$, it holds
\begin{equation}
\label{liminf}
G(x)\le\liminf_{n\to\infty}G_n(x_n).
\end{equation}
\vspace*{-9mm}
\item[(ii)] (Recovery sequence)\\
For every $x\in X$, there exists a sequence $(x_n)_{n\in\N}\subset X$ with
$x_n\to x$ in $X$ such that
\[ G(x)\ge\limsup_{n\to\infty}G_n(x_n). \]
\end{itemize}
\end{Definition}

\begin{Definition}
\label{def2} Let $X$ be a topological space.
\vspace*{-2mm}
\begin{itemize}
\item[(i)] A functional $G:X\to\overline{\R}$ is {\em coercive} on $X$ if for
all $\alpha\in\R$ the closure of the sublevel sets
$\{x\in X\;|\; G(x)\le\alpha\}$ is compact in X.
\vspace*{-2mm}
\item[(ii)] A family of functionals $\{G_n\}_{n\in\N}:X\to\overline{\R}$ is
{\em equi-coercive} if for all $\alpha\in\R$ there exists a compact set
$K_\alpha\subset X$ such that for all $n\in\N$,
$\{x\in X\;|\;G_n(x)\le\alpha\}\subset K_\alpha$.
\end{itemize}
\end{Definition}
If $G$ is coercive there exists a compact set $K\subset X$ with
$\inf_{x\in X}G(x)=\inf_{x\in K}G(x)$. If $(G_n)_{n\in\N}$ is equi-coercive
there exists a compact set $K\subset X$ such that for all $n\in\N$
$\inf_{x\in X}G_n(x)=\inf_{x\in K}G_n(x)$. A family $(G_n)_{n\in\N}$ is
equi-coercive iff there exists a lower semicontinuous and coercive
$\psi:X\to\overline{\R}$ such that $G_n\ge\psi$ for all $n\in\N$.

The choice of topology in $X$ is crucial: The finer the topology in $X$, the
easier the lower semicontinuity property inherent in (\ref{liminf}) is
satisfied. Contrary, the compactness of the sublevel sets $K_\alpha$ calls for
coarser topologies. Hence, both conditions are competing. Often, the weak
topology in a Sobolev space provides a good compromise where both properties
are satisfied simultaneously.

\vspace*{2mm}
\begin{Theorem}[Fundamental theorem of Gamma-convergence]
\label{theo1}
Let $X$ be a topological space and $(G_n)_{n\in\N}:X\to\overline{\R}$ be an
equi-coercive family of functionals with $G_n\to G$ in the $\Gamma$-sense.
Then it holds
\begin{itemize}
\item[(i)] $G$ is coercive.
\vspace*{-2mm}
\item[(ii)] The minima of $G_n$ converge to the minima of $G$, i.e. letting
$\underline{e}_n:=\inf_{x\in X}G_n(x)$, $\underline{e}:=\inf_{x\in X}G(x)$, then
\[ \underline{e}_n\to\underline{e}\qquad\mbox{for }n\to\infty. \]
\item[(iii)] The minimizers of $G_n$ converge to the minimizers of $G$.
\end{itemize}
\end{Theorem}

\vspace*{2mm}
\noindent{\bf Proof.} See, e.g., \cite{DalMaso93}. \qed

\section{Zero-order Gamma-limit}
\label{seczero}
We derive the first two terms in the Gamma-expansion of $E_\ve$ beginning
in this section with the zero order $\Gamma$-limit. Writing $z:=u'$, let
(cf. Eqn.~(\ref{Eeps}))
\begin{align}
Q(z,\alpha) \;:=\;& \frac{\mu}{2}|z|^2+W(z,\alpha)\nn\\
\label{Qdef}
=\;& \frac{\mu}{2}|z|^2+\frac{\mu}{2}\Big(\sin(\alpha)z-4\sin^2\!\big(
\frac{\alpha}{2}\big)\Big)^2+\frac{\mu_c}{2}\Big(\cos(\alpha)z-2\sin(\alpha)
\Big)^2.
\end{align}

\begin{figure}[h!tb]
\unitlength1cm
\begin{picture}(11.5,4.8)
\put(4.5,-0.3){\psfig{figure=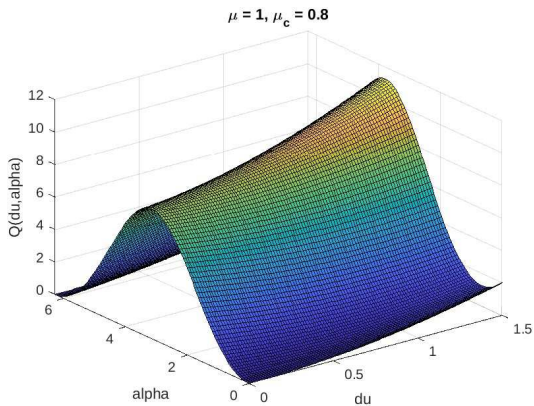,width=7.0cm}}
\end{picture}
\caption{\label{fig3}
Plot of $Q(z,\alpha)$ for $\alpha\in[0,2\pi]$, $z\in[0,1.5]$,
$\mu=1$ and $\mu_c=0.8$.}
\end{figure}

\begin{Definition}
\label{def3}
For $Q$ given by (\ref{Qdef}) we denote by $Q^{**}$ the convex envelope of
$Q$, i.e.
\begin{equation}
\label{gss}
Q^{**}(z,\alpha):=\sup\big\{g(z,\alpha)\;|\;g
\mbox{ is convex and }g(z,\alpha)\le Q(z,\alpha)\big\}.
\end{equation}
\end{Definition}
Since $Q$ is finite it is not necessary to demand in (\ref{gss}) that $g$ be
lower semi-continuous as finite convex functions are automatically continuous.

\vspace*{2mm}
The computation of $Q^{**}$ is postponed to Lemma~\ref{lem2} below.

\vspace*{2mm}
We define the relaxed functional
\begin{equation}
\label{E0def}
E_0(u,\alpha):=\left\{\!\!\begin{array}{ll}
\int_0^1Q^{**}(u',\alpha)\dx, &
\mbox{ if }(u,\alpha)\in {\cal X},\\
+\infty, & \mbox{ else.} \end{array}\right.
\end{equation}

In order to single out a solution, for given $\theta\in[0,2\pi]$, we
introduce the volume constraint
\begin{equation}
\label{vc}
\int_0^1\alpha(x)\dx=\theta
\end{equation}
and the corresponding functional for $\ve>0$
\begin{equation}
\label{Endef}
E_{\ve}^\theta(u,\alpha):=\left\{\!\!\!\begin{array}{ll}
\int_0^1\ve^2|\alpha'|^2\!+\!Q(u',\alpha)\dx, &
\mbox{if }(u,\alpha)\in{\cal X},\,\int_0^1\alpha(x)\dx=\theta,\\
+\infty, & \mbox{else.} \end{array}\right.
\end{equation}

\newpage
\begin{Lemma}
\label{lem1}
Let $(\ve_n)_{n\in\N}$ be a sequence of positive real numbers with
$\ve_{n}\searrow0$ for $n\to\infty$. Then the family of functionals
$(E_{\ve_n}^\theta)_{n\in\N}$ defined in (\ref{Endef}) is equi-coercive with
respect to weak convergence in $(L^2(\Omega))^2$.
\end{Lemma}

\noindent{\bf Proof.}
We need to show that for every sequence
$(u_n,\alpha_n)_{n\in\N}\subset(L^2(\Omega))^2$ with
\linebreak
$\sup_{n\in\N}E_{\ve_n}^\theta(u_n,\alpha_n)<\infty$ there exists a subsequence
$(u_{n_k},\alpha_{n_k})_{k\in\N}$ that converges weakly in $(L^2(\Omega))^2$.
Directly from (\ref{Endef}) and (\ref{Qdef}) we obtain 
\begin{equation}
\label{est}
\frac{\mu}{2}\|u_n'\|_{L^2(\Omega)}^2
\le E_{\ve_n}^\theta(u_n,\alpha_n)\le C
\end{equation}
uniformly in $n\in\N$. From the definition (\ref{Endef}) of $E_{\ve_n}^\theta$,
this implies
$(u_n,\alpha_n)\in{\cal X}={\cal X}_u\times{\cal X}_\alpha$ for every $n\in\N$.
Functions $u\in{\cal X}_u$ satisfy for $\uhom(x):=\gamma\,x$ and $C_0=0$
the cone condition:
\[ \mbox{If }\uhom+\xi\in{\cal X}_u \mbox{ for }\xi\in\R, \mbox{ then }
|\xi|\le C_0. \]
Hence, the general Poincar{\'e} inequality applies for functions in
${\cal X}_u$ and from (\ref{est}) and $\mu>0$, we infer that
$\|u_n\|_{W^{1,2}(\Omega)}$ is bounded uniformly in $n$.
By the Banach-Alaoglu theorem, there exists a subsequence $(u_{n_k})_{k\in\N}$
that converges weakly in $L^2(\Omega)$.

From $(\alpha_n)_{n\in\N}\subset{\cal X}_\alpha$ follows
$\alpha_n\in[0,2\pi]$ pointwise in $\Omega$ and hence the boundedness
$\|\alpha_n\|_{L^2(\Omega)}\le2\pi$ uniformly in $n\in\N$.
Again, (see, e.g., Theorem~2.6 in \cite{Braides02})
this implies the existence of a subsequence $(\alpha_{n_k})$ converging
to $\alpha$ weakly in $L^2(\Omega)$.
Finally, the constraint (\ref{vc}) is closed with respect to weak convergence in
$L^2(\Omega)$ as $\int_0^1\big(\alpha_n\!-\!\alpha)\zeta\dx\to0$ for
$n\to\infty$ holds for any test function $\zeta\in L^2(\Omega)$.
So the limit function $\alpha$ satisfies (\ref{vc}). \qed

\vspace*{4mm}
As a consequence of Lemma~\ref{lem1}, all the statements of Theorem~\ref{theo1}
hold for $(E_{\ve_n})_{n\in\N}$ and the limit functional $E_0$.

\newpage
\begin{Proposition}[Zero-order Gamma-limit]
\label{prop1}
Let $(\ve_n)_{n\in\N}$ be a sequence of positive real numbers converging to $0$.
For $\theta\!\in\![0,2\pi]$, let $E_{\ve_n}^\theta$ on $(L^2(\Omega))^2$ be
given by (\ref{Endef}),
\[ E_0(u,\alpha)=\left\{\!\!\begin{array}{ll} \int_0^1Q^{**}(u',\alpha)\dx, &
\mbox{ if }(u,\alpha)\in {\cal X},\\
+\infty, & \mbox{ else.} \end{array}\right. \]
Then the $\Gamma$-limit of $E_{\ve_n}^\theta$ with respect to weak convergence
in $L^2(\Omega)$ for $n\to\infty$ is
\begin{equation}
\label{E0thetadef}
E_0^\theta(u,\alpha)=\left\{\begin{array}{ll}E_0(u,\alpha), &
\mbox{ if }(u,\alpha)\in{\cal X},\quad\int_0^1\alpha(x)\dx=\theta,\\
+\infty, & \mbox{ else.} \end{array}\right.
\end{equation}
\end{Proposition}
\noindent{\bf Proof.}
(i) {\it Liminf inequality}.
We may restrict our attention to a sequence
$(u_n,\alpha_n)_{n\in\N}\subset{\cal X}$
converging to $(u,\alpha)\in{\cal X} $ as otherwise
$\liminf_{n\to\infty}E_{\ve_n}(u_n,\alpha_n)=+\infty$. Since the constraint
$\int_0^1\alpha(x)\dx\!=\!\theta$ is closed w.r.t. weak convergence in
$L^2(\Omega)$, we find
\vspace*{-2mm}
\begin{eqnarray}
\liminf_{n\to\infty}E_{\ve_n}(u_n,\alpha_n) &=& \liminf_{n\to\infty}
\int_0^1\ve_n^2|\alpha_n'|^2+Q(u_n',\alpha_n)\dx\nn\\
\label{lowinf}
&\ge& \liminf_{n\to\infty}
\int_0^1\ve_n^2|\alpha_n'|^2+Q^{**}(u_n',\alpha_n)\dx.
\end{eqnarray}
Here, (\ref{gss}) was used to get to the second line.

From its definition in (\ref{gss}), $Q^{**}$ is convex and lower semicontinuous
w.r.t. strong convergence. Hence its epigraph is convex and closed, hence
weakly closed. This shows that $Q^{**}$ is weakly lower semicontinuous.
This yields the desired inequality
\[ \liminf_{n\to\infty}\int_0^1\ve_n|\alpha_n'|^2+Q^{**}(u_n',\alpha_n)\dx
\;\ge\; \int_0^1Q^{**}(u',\alpha)\dx=E_0(u,\alpha). \]

\vspace*{2mm}
\noindent (ii) {\it Limsup inequality}.

If $(u,\alpha)\notin{\cal X}$, then $E_0^\theta(u,\alpha)=E_0(u,\alpha)=+\infty$
and the Limsup inequality is obviously satisfied.
For $(u,\alpha)\in{\cal X}$, we may simply choose $\alpha_n:=\alpha$,
$u_n:=u$ for $n\in\N$ as a recovery sequence such that
\begin{equation}
\label{ls1}
\Gamma\!-\!\limsup_{n\to\infty}E_{\ve_n}(u,\alpha)\le\int_0^1
Q(u',\alpha)\dx\qquad\mbox{in }{\cal X}.
\end{equation}
By density,
c.f. \cite[Remark~1.29]{Braides02},
(\ref{ls1}) remains valid in $(L^2(\Omega))^2$. Taking the convex
lower semi-continuous hull on both sides of (\ref{ls1}) immediately yields 
the limsup inequality
\[ \limsup_{n\to\infty}E_{\ve_n}(u_n,\alpha_n)=\limsup_{n\to\infty}
E_{\ve_n}(u,\alpha)\le\int_0^1Q^{**}(u',\alpha)\dx=E_0(u,\alpha). \]

\noindent From $E_0^\theta(u,\alpha)<\infty$ we see that each $\alpha_n=\alpha$
satisfies the volume constraint (\ref{vc}). \qed

\newpage
For the computation of $Q^{**}$ and the first-order Gamma-limit,
we need to consider the minimal energy $\emin$ of $E$ (defined in (\ref{Edef})),
\[ \emin:=\min\Big\{E(u,\alpha)\;\Big|\;(u,\alpha)\in{\cal X}\Big\}. \]
The following Table~\ref{tab1} has the details.

{\renewcommand{\arraystretch}{1.05}
\begin{table}[h!bt]
\hspace*{10pt}
\begin{tabular}{|l|l|l|}\hline
\multicolumn{1}{|c|}{\mbox{\footnotesize Parameters}}&
\multicolumn{1}{|c|}{\mbox{\footnotesize Minimizers}}&
\multicolumn{1}{|c|}{\mbox{\footnotesize Minimal energy $\emin$ of $E$}}\\
\hline\hline
$\mu=\mu_c$ & $(\uhom,\alpha_2)$ &
$\mu\big(\gamma^2+4-2\sqrt{\gamma^2+4}\big)$\\
$\mu_c=0$   & $(\uhom,\AA\!=\!0),\,(\uhom,\BB\!=\!\arctan(\!
\frac{4\gamma}{4-\gamma^2}\!))$ &$\frac{\mu}{2}\gamma^2$\\
$\mu_c>\muccrit$, $\mu\not=\mu_c$ & $(\uhom,\alpha_2)$ &
$\mu\big(\gamma^2+4-2\sqrt{\gamma^2+4}\big)$\\
$\mu>\mu_c$, $0<\mu_c\le\muccrit$ & $(\uhom,\AA)$, $(\uhom,\BB)$ &
$\frac{\mu+\mu_c}{2}\gamma^2-\frac{2\mu_c^2}{\mu-\mu_c}$\\
\hline
\end{tabular}
\caption{\label{tab1}
Unique minimizers of $E$ for $L_c=0$ in ${\cal X}$ and minimal energies
for different parameter ranges.}
\end{table}}

The data of Table~\ref{tab1} is taken from \cite{BN22}. We adopt the notations
\begin{equation}
\label{uhomdef}
\uhom(x):=\gamma\,x,\qquad 0\le x\le 1
\end{equation}
for the {\it homogeneous deformation} which turns out to be optimal in
${\cal X}_u$, and
\begin{equation}
\label{ABdef}
\AA:=\arctan\Big(\frac{\gamma\mu-f}{2\mu+\frac{\gamma}{2}f}\Big),\qquad
\BB:=\arctan\Big(\frac{\gamma\mu+f}{2\mu-\frac{\gamma}{2}f}\Big),\qquad
\alpha_2:=\arctan\!\Big(\frac{\gamma}{2}\Big)
\end{equation}
for the global minimizers of $E$ with $L_c=0$, where
\begin{equation}
\label{fdef}
f:=\Big((\gamma^2+4)(\mu-\mu_c)^2-4\mu^2\Big)^{1/2}.
\end{equation}
The {\it critical value of} $\mu_c$ is given by
\begin{equation}
\label{muccritdef}
\muccrit:=\mu\Big[1-\frac{2}{\sqrt{\gamma^2+4}}\Big].
\end{equation}

\begin{Remark}
\label{rem1}
The case $\mu_c=0$ is a limiting case of the regime $\mu>\mu_c$,
$0<\mu_c\le\muccrit$. For $\mu_c=0$, by Eqn.~(\ref{fdef}), $f=\mu\gamma$ such
that by Formula~(\ref{ABdef})
\begin{equation}
\label{alphadef}
\AA=0,\qquad \BB=\arctan\!\Big(\frac{4\gamma}{4-\gamma^2}\Big).
\end{equation}
This formula had already been derived in \cite{BN22}.
Therein, the value $\BB$ for $\mu_c=0$ had been denoted $\alpha_3$ and
introduced by the identity
\begin{equation}
\label{a3def}
\alpha_3=\eta^{-1}(\gamma),\qquad \eta(\alpha):=\frac{4\sin^2\!\big(
\frac{\alpha}{2}\big)}{\sin(\alpha)}.
\end{equation}
The inverse $\eta^{-1}(\gamma)$ in (\ref{a3def}) exists for $0\le\gamma<2\pi$
while $\BB$ in (\ref{alphadef}) exists for $\gamma\in[0,2)$.
\end{Remark}

\begin{figure}[h!p]
\unitlength1cm
\begin{picture}(11.5,10.8)
\put(0.5,-0.2){\psfig{figure=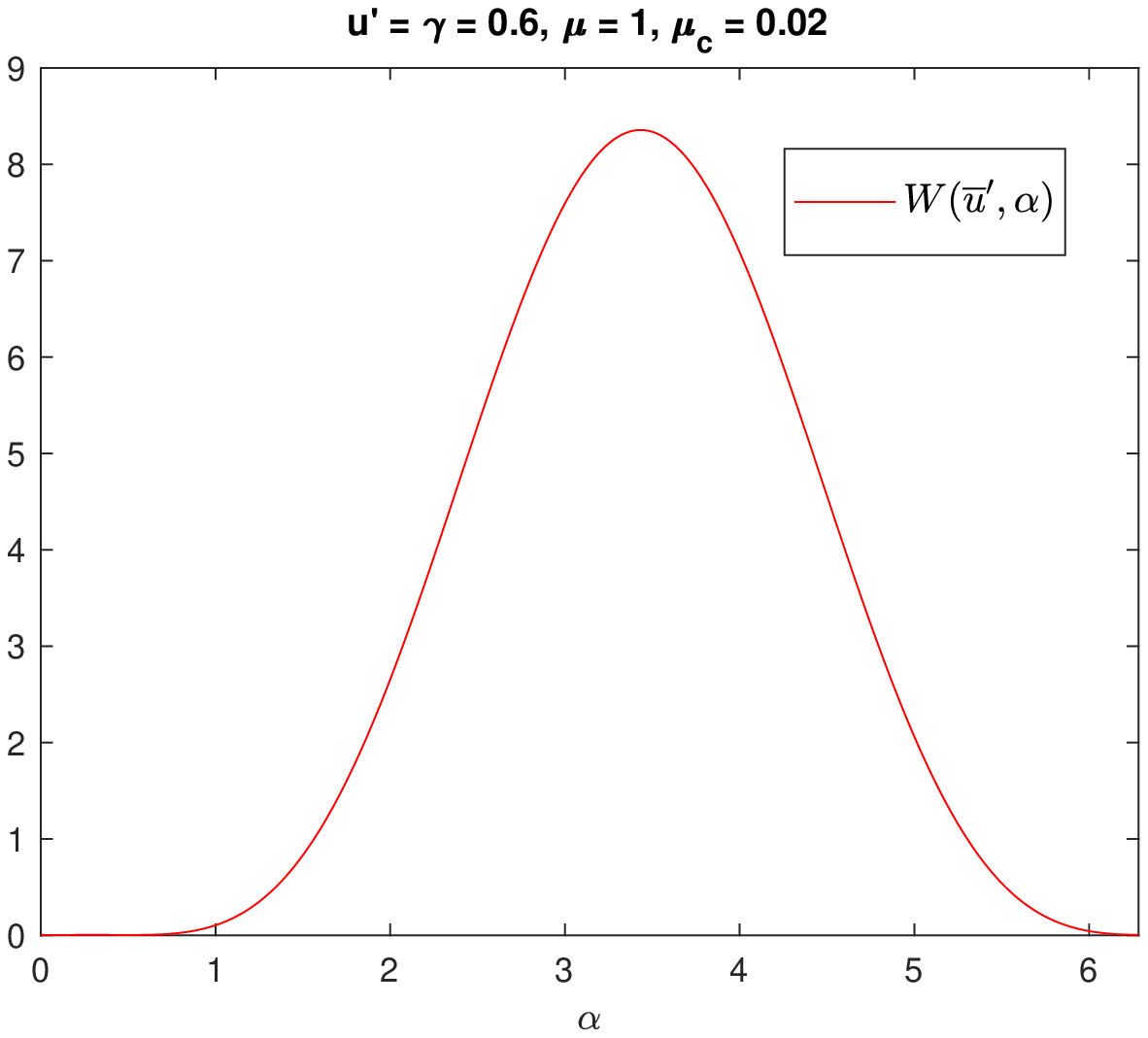,width=7.0cm}}
\put(8.0,-0.2){\psfig{figure=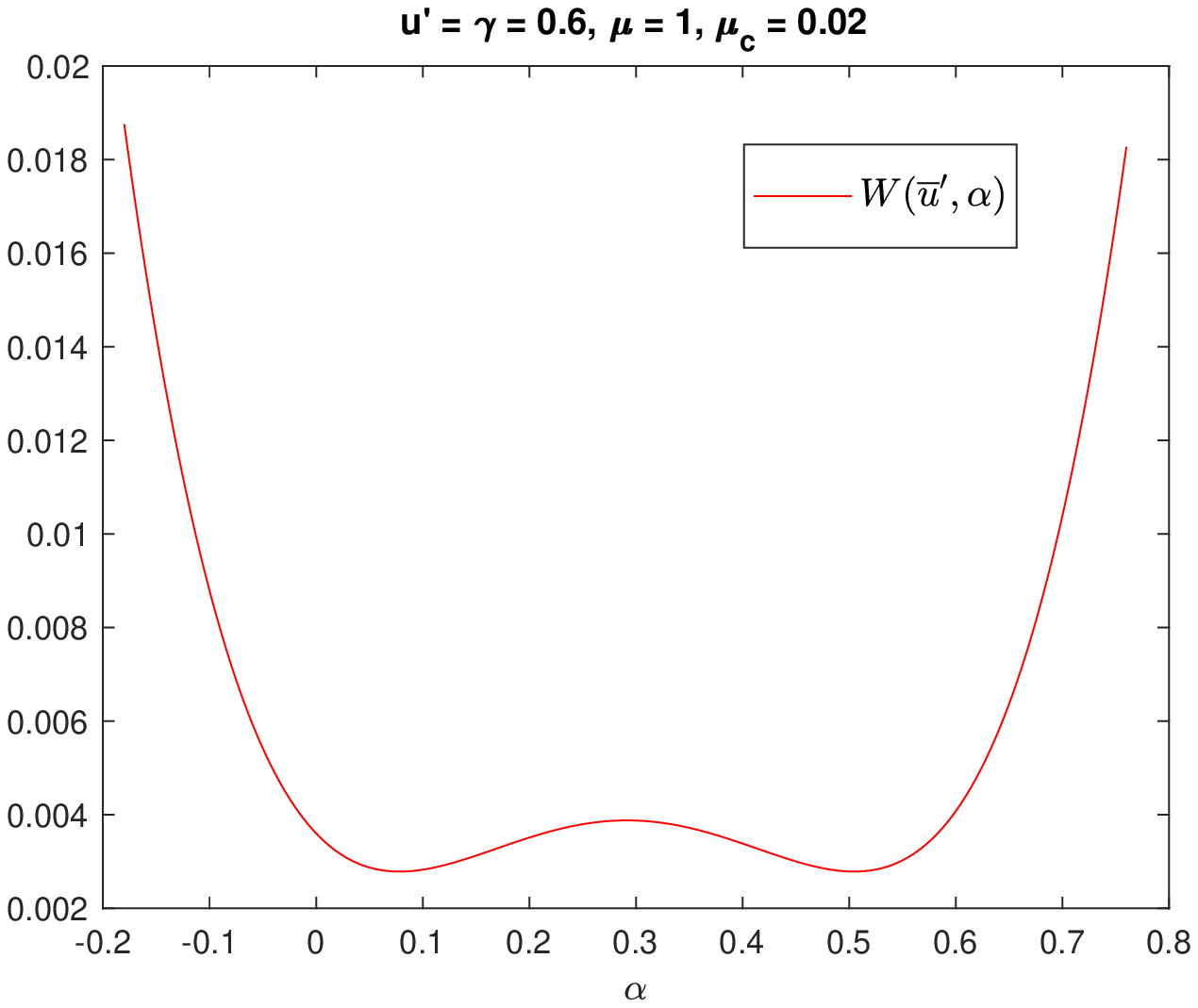,width=7.0cm}}
\put(0.5,5.3){\psfig{figure=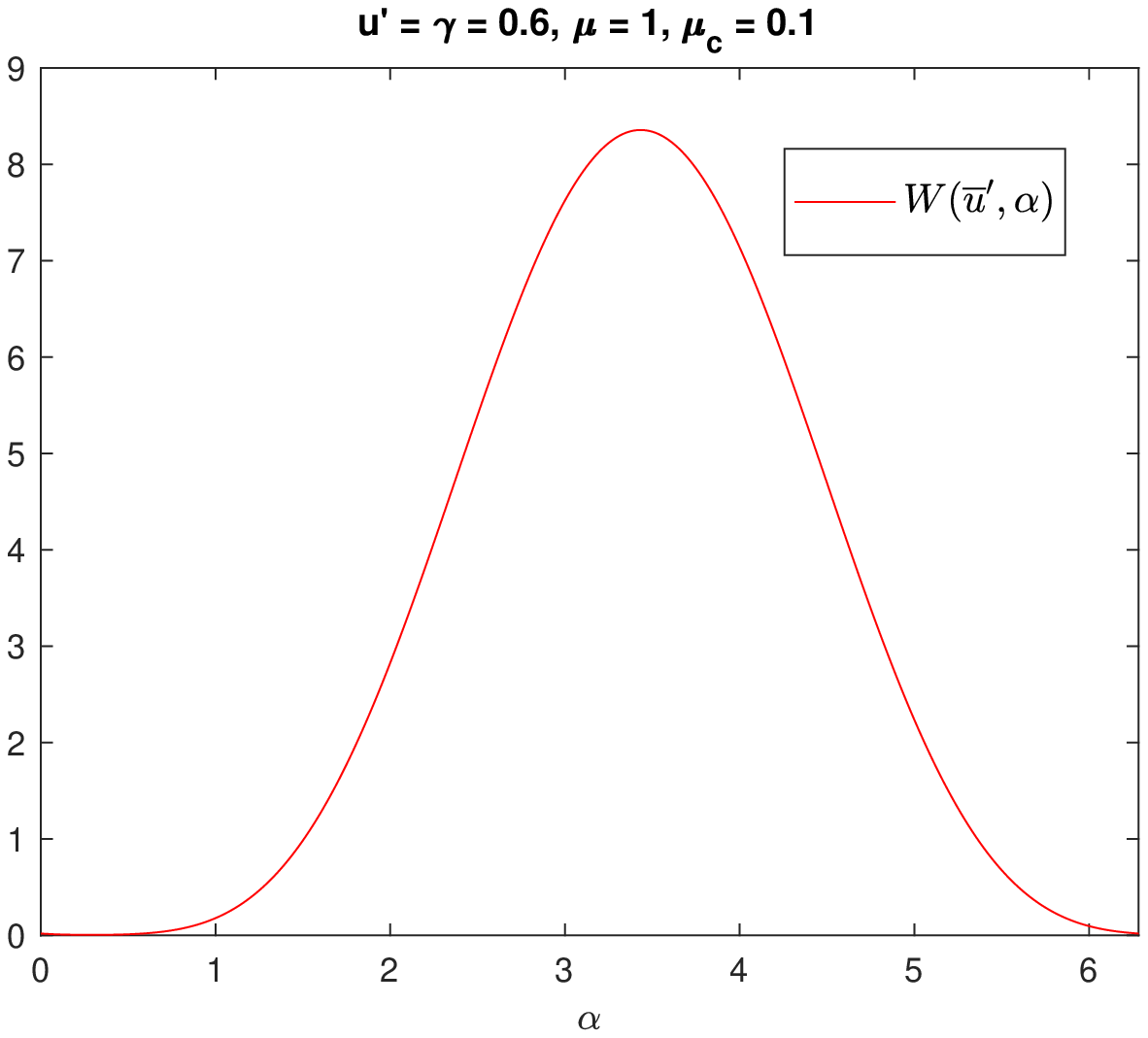,width=7.0cm}}
\put(8.0,5.3){\psfig{figure=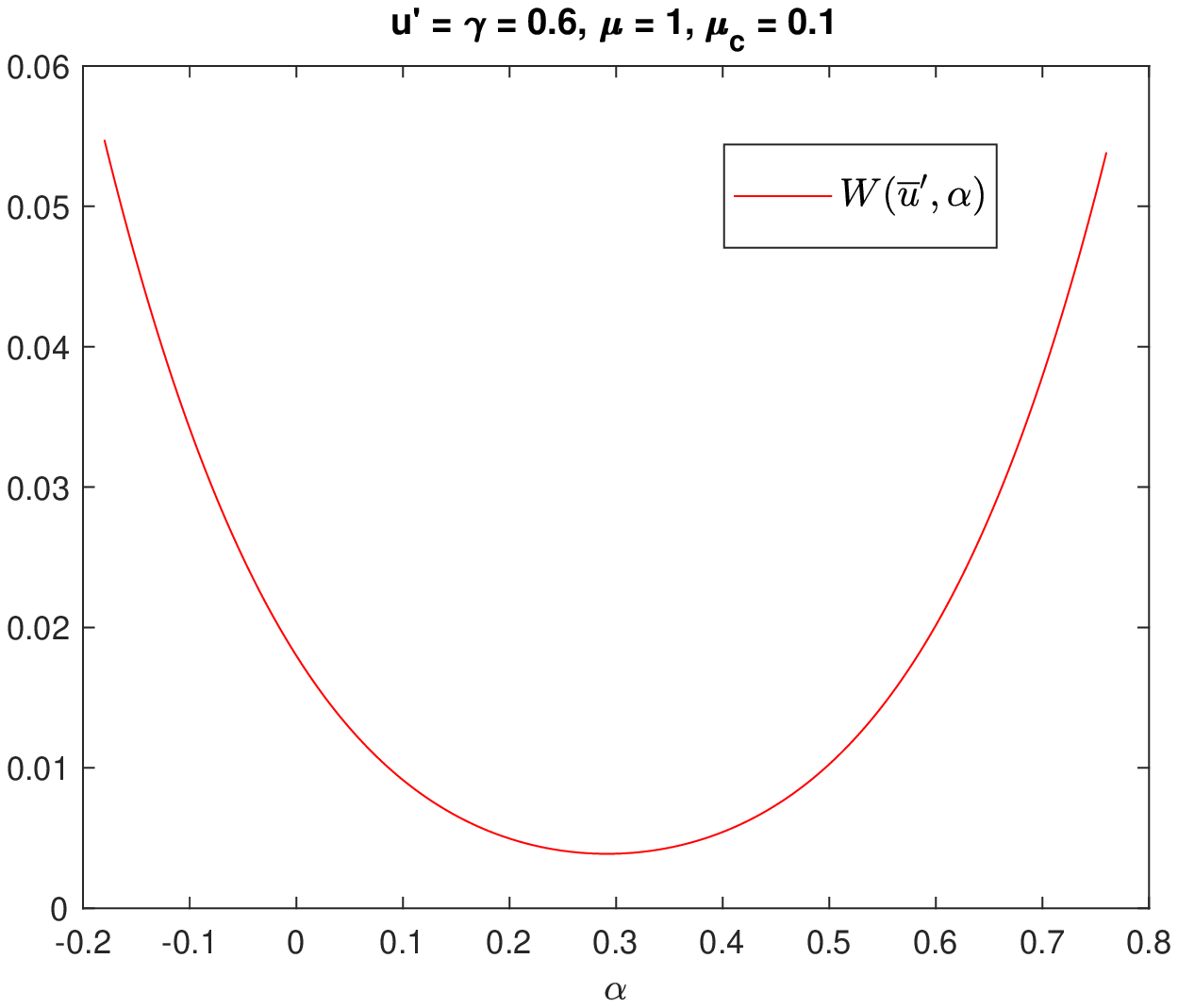,width=7.0cm}}
\put(10.10,0.0){{\tiny $\alpha_1^{\!-}$}}
\put(10.34,0.31){\line(0,1){0.15}}
\put(12.69,0.0){{\tiny $\BB$}}
\put(12.87,0.31){\line(0,1){0.15}}
\put(11.25,5.54){{\tiny $\alpha_2$}}
\put(11.54,5.80){\line(0,1){0.15}}
\put(7.945,2.70){\color{blue}\line(-1,1){0.2}}
\put(7.1,2.7){\color{blue}\line(2,0){0.85}}
\put(7.945,2.70){\color{blue}\line(-1,-1){0.3}}
\put(1.2,0.3){\color{blue}\line(0,2){0.3}}
\put(2.1,0.3){\color{blue}\line(0,2){0.3}}
\put(1.2,0.3){\color{blue}\line(2,0){0.9}}
\put(1.2,0.6){\color{blue}\line(2,0){0.9}}
\put(8.2,-0.2){\color{blue}\line(0,2){5.4}}
\put(14.6,-0.2){\color{blue}\line(0,2){5.4}}
\put(8.2,-0.2){\color{blue}\line(2,0){6.4}}
\put(8.2,5.2){\color{blue}\line(2,0){6.4}}
\end{picture}
\caption{\label{fig4}
\footnotesize Plots of $\alpha\mapsto W(\uhom',\alpha)$ for
$\mu=1$ and $\gamma=0.6$.
Left: Plots for $\alpha\in[0,2\pi]$. Right: Close-ups for $\alpha\in[-0.2,0.8]$.
Bottom: For $\mu_c=0.02<\muccrit=0.0422$ there is a double-well with minimal
value $\wmin=\frac{\mu_c}{2}\gamma^2-\frac{2\mu_c^2}{\mu-\mu_c}\approx
0.00278$ at $\AA\approx0.0783$ and at $\BB\approx0.5046$ in accordance with
Eqn.~(\ref{ABdef}).
Top: Strict convexity in the interval $\alpha\in[0,1]$ for
$\mu_c=0.1>\muccrit$ with the unique minimal value
$\wmin=\mu\big[\frac{\gamma^2}{2}+4-2\sqrt{\gamma^2+4}\big]\approx0.003877$ at
$\alpha_2=\arctan(\gamma/2)\approx0.2915$. The double well is very flat and
practically invisible in the full plot displayed on the bottom line left.
The blue box on the bottom left illustrates the section which is enlarged
on the right.}
\end{figure}

\begin{figure}[h!tb]
\unitlength1cm
\begin{picture}(11.5,4.7)
\put(0.5,-0.1){\psfig{figure=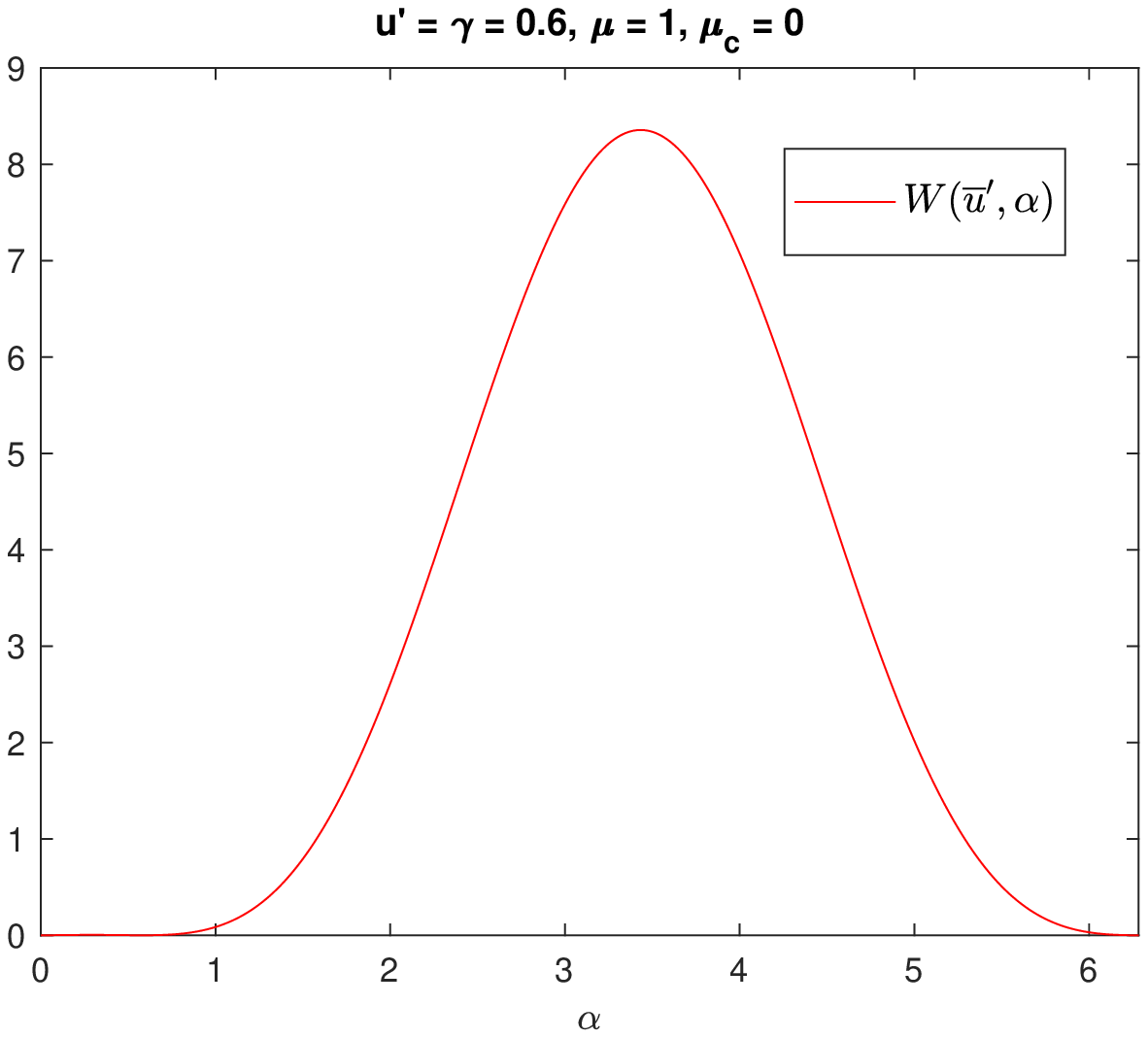,width=7.0cm}}
\put(8.0,-0.1){\psfig{figure=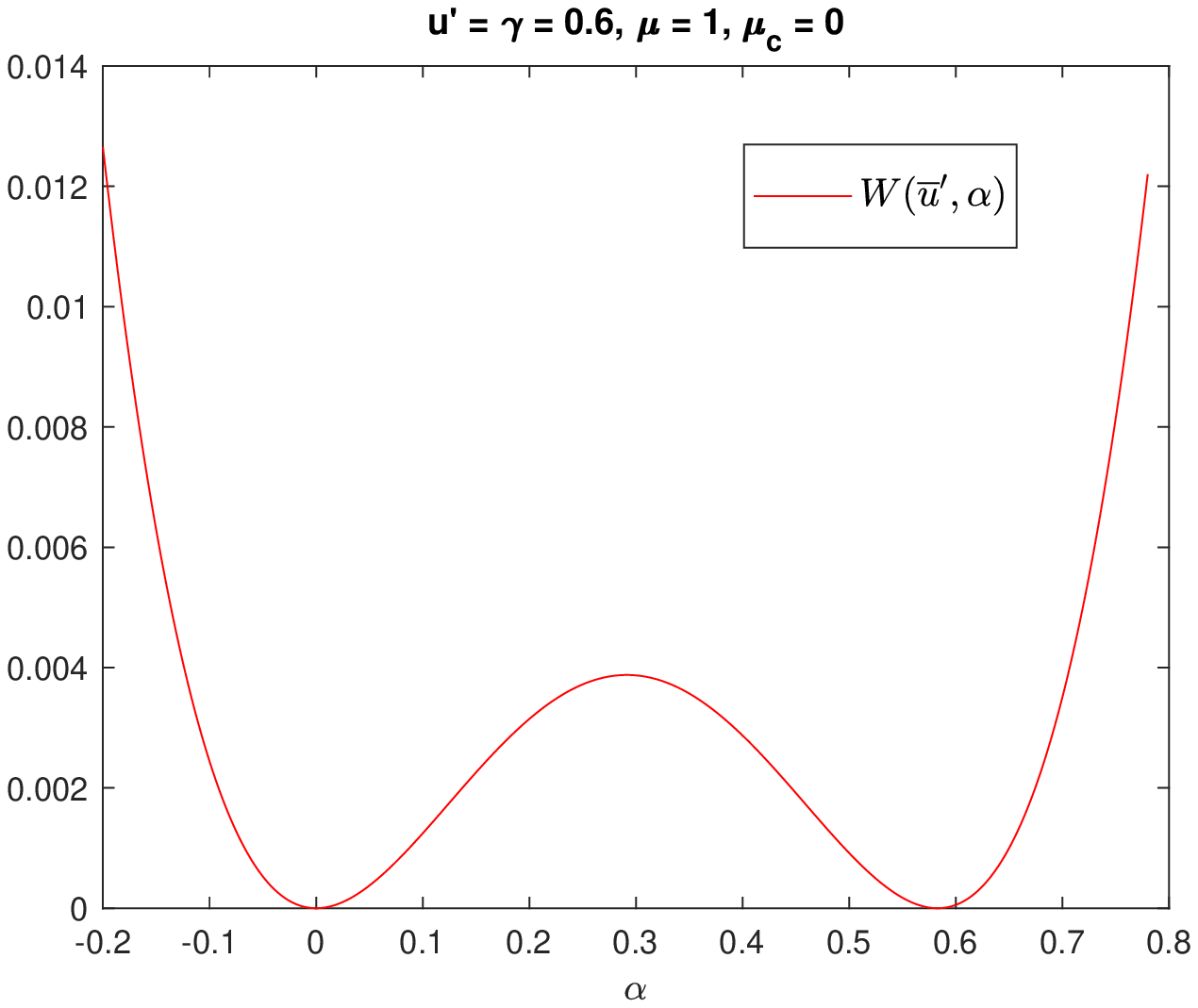,width=7.0cm}}
\put(1.2,0.3){\color{blue}\line(0,2){0.3}}
\put(2.1,0.3){\color{blue}\line(0,2){0.3}}
\put(1.2,0.3){\color{blue}\line(2,0){0.9}}
\put(1.2,0.6){\color{blue}\line(2,0){0.9}}
\put(7.945,2.70){\color{blue}\line(-1,1){0.2}}
\put(7.1,2.7){\color{blue}\line(2,0){0.85}}
\put(7.945,2.70){\color{blue}\line(-1,-1){0.3}}
\put(8.2,-0.2){\color{blue}\line(0,2){5.4}}
\put(14.6,-0.2){\color{blue}\line(0,2){5.4}}
\put(8.2,-0.2){\color{blue}\line(2,0){6.4}}
\put(8.2,5.2){\color{blue}\line(2,0){6.4}}
\put(9.85,0.04){{\tiny $\AA$}}
\put(9.99,0.40){\line(0,1){0.15}}
\put(13.00,0.04){{\tiny $\BB$}}
\put(13.15,0.40){\line(0,1){0.15}}
\end{picture}
\caption{\label{fig5}
\footnotesize Two plots of $\alpha\mapsto W(\uhom',\alpha)$ for $\mu_c=0$,
$\mu=1$ and $\gamma=0.6$. Left: Plot for $\alpha\!\in\![0,2\pi]$.
Right: Close-up for $\alpha\in[-0.2,0.8]$ showing a double-well with minima
at $\AA=0$ and
$\BB=\arctan\!\big(\frac{4\gamma}{4-\gamma^2}\!\big)\approx0.5829$
as predicted by Eqn.~(\ref{alphadef}). The double well is very flat and does
not show up in the full plot on the left.}
\end{figure}

\newpage
In the interval $\alpha\in[0,1]$ which contains all minimizers,
$W(\uhom',\cdot)$ is strictly convex for $\mu_c>\muccrit$,
while for $\mu_c=0$ and in the non-classical regime $\mu>\mu_c$ with
$0<\mu_c\le\muccrit$, $W(\uhom',\cdot)$
is a double-well potential with minimizers at $\AA$ and at $\BB$,
cf. Fig.~\ref{fig4}. The limiting case $\mu_c=0$ is displayed in
Fig.~\ref{fig5}.

\newpage
To complete the discussion of the zero-order $\Gamma$-limit, it remains to
compute the convexification $Q^{**}$.
\begin{Lemma}[Computation of $Q^{**}$]
\label{lem2}
Let $Q$ be given by (\ref{Qdef}). For $z>0$, depending on $\mu$ and $\mu_c$,
the convexification $Q^{**}(z,\alpha)$ is given by the following formulas.

\vspace*{2mm}
\noindent(i) If $\mu=\mu_c$ or $\mu_c>\muccrit:$
\begin{equation}
\label{Qss1}
Q^{**}(z,\alpha) \;=\; \left\{\begin{array}{ll}
W(z,\alpha)+\frac{\mu}{2}z^2 & \mbox{ if }\alpha\in[0,\alpha_2),\\
\big(\frac{\mu+\mu_c}{2}z^2-\emin\big)\frac{\alpha-\alpha_2}{2\pi-\alpha_2}
+\emin & \mbox{ if }\alpha\in[\alpha_2,2\pi) \end{array}\right.
\end{equation}
where $\emin=\emin(z):=\mu\big(z^2+4-2\sqrt{z^2+4}\big)$ and
$\alpha_2=\alpha_2(z):=\arctan\!\big(\frac{z}{2}\big)$, cf. Eqn.~(\ref{ABdef}).

\vspace*{2mm}
\noindent(ii) If $\mu_c=0:$
\begin{equation}
\label{Qss2}
Q^{**}(z,\alpha)=\frac{\mu}{2}z^2.
\end{equation}

\vspace*{2mm}
\noindent(iii) If $0<\mu_c\le\muccrit$, $\mu\not=\mu_c:$
\begin{equation}
\label{Qss3}
Q^{**}(z,\alpha) \;=\; \left\{\begin{array}{ll}
W(z,\alpha)+\frac{\mu}{2}z^2 & \mbox{ if }\alpha\in[0,\AA),\\
\emin & \mbox{ if }\alpha\in[\AA,\BB),\\
\big(\frac{\mu+\mu_c}{2}z^2-\emin\big)\frac{\alpha-\BB}{2\pi-\BB}+\emin
& \mbox{ if }\alpha\in[\BB,2\pi), \end{array}\right.
\end{equation}
\vspace*{-2mm}
where
$\emin\!=\!\emin(z)\!=\!\frac{\mu+\mu_c}{2}z^2\!-\!\frac{2\mu_c^2}{\mu-\mu_c}$,
$\AA\!=\!\AA(z)\!=\!\arctan\!\big(\!\frac{\gamma\mu-f}{2\mu+\frac{z}{2}f}
\!\big)$, $\BB\!=\!\BB(z)\!=\!\arctan\!\big(\!
\frac{\gamma\mu+f}{2\mu-\frac{z}{2}f}\!\big)$ for
$f:=\big((z^2+4)(\mu-\mu_c)^2-4\mu^2\big)^{1/2}$, cf. Eqn.~(\ref{ABdef}).
\end{Lemma}

\vspace*{2mm}
\noindent{\bf Proof}.
The convex envelope can be computed by the formula
\begin{align*}
Q^{**}(z,\alpha) \;=\; \min\Big\{tQ(z_1,\alpha_1)+(1-t)Q(z_2,\alpha_2)\;\Big|\;&
tz_1+(1-t)z_2=z,\\
& t\alpha_1+(1-t)\alpha_2=\alpha,\;t\in[0,1]\Big\}.
\end{align*}
From $Q(z,\alpha)=\frac\mu2|z|^2+W(z,\alpha)$, 
the definition of $Q^{**}$ in (\ref{gss}),
and the convexity of $\frac{\mu}{2}|z|^2$, we have
\[ Q^{**}(z,\alpha)=\frac{\mu}{2}|z|^2+W^{**}(z,\alpha). \]
Furthermore,
\[ \frac{\partial^2 Q}{\partial z^2}(z,\alpha)\;=\; \mu\big[1+\sin^2(\alpha)
\big]+\mu_c\cos^2(\alpha)>0. \]
Hence the function $z\mapsto Q(z,\alpha)$ is strictly monotonically increasing
and convex.

A function is convex iff its restriction to any line that intersects its domain
is convex. Similarly, a function $g$ in Definition~(\ref{gss}) is convex iff
its two coordinate functions $z\mapsto g(z,\alpha)$ for fixed $\alpha$ and
$\alpha\mapsto g(z,\alpha)$ for fixed $z$ are convex, see, e.g.,
\cite[p.~67]{Boyd09}.
Therefore it only remains to compute the convexification of
$\alpha\mapsto W(z,\alpha)$ for fixed $z$, see Fig.~\ref{fig3} for an
illustration of $Q(z,\alpha)$.

The minimizers of $W(z,\cdot)$ for fixed $z$ are readily available by
Table~\ref{tab1} and Eqns.~(\ref{ABdef}), (\ref{fdef}).
Therein, we only need to replace $\gamma$ by $z$.
For the rest of the proof, as $z$ is fixed, we write for short
$\alpha_1^\pm$, $\alpha_2$ and $\emin$ instead of
$\alpha_1^\pm(z)$, $\alpha_2(z)$ and $\emin(z)$.

\begin{figure}[phtb]
\unitlength1cm
\begin{picture}(11.5,15.0)
\put(2.0,7.6){\psfig{figure=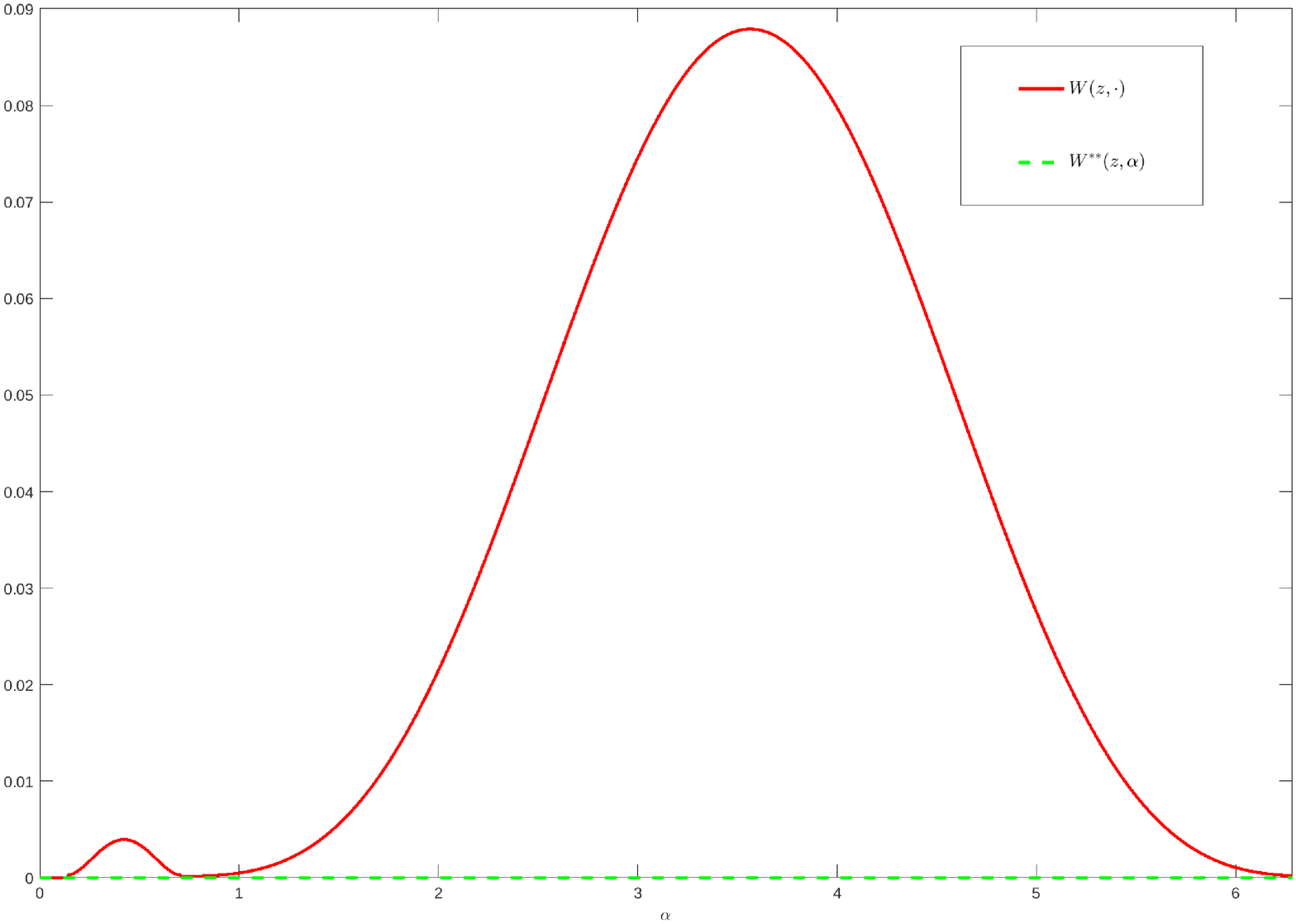,width=13.2cm}}
\put(2.0,-0.3){\psfig{figure=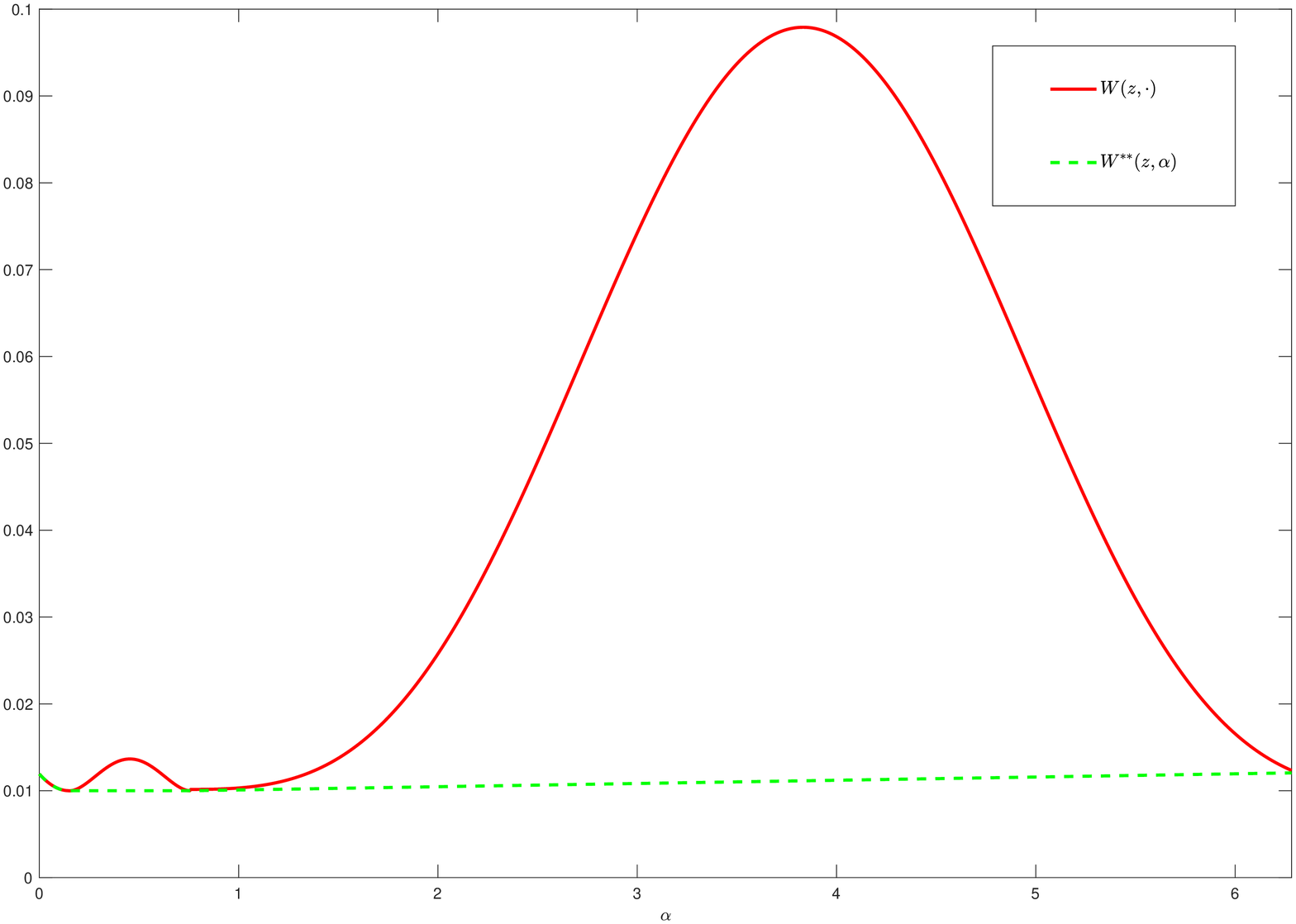,width=13.2cm}}
\put(4.00,0.56){\line(0,1){0.15}}
\put(4.95,0.56){\line(0,1){0.15}}
\put(3.80,0.25){{\tiny $\AA$}}
\put(4.80,0.25){{\tiny $\BB$}}
\put(13.80,0.25){{\tiny $2\pi$}}
\put(3.70,8.45){\line(0,1){0.15}}
\put(4.90,8.45){\line(0,1){0.15}}
\put(3.50,8.15){{\tiny $\AA$}}
\put(4.80,8.15){{\tiny $\BB$}}
\put(13.80,8.15){{\tiny $2\pi$}}
\end{picture}
\caption{\label{fig6}
\footnotesize Sketch of the construction of $W^{**}(z,\cdot)$.
Top: The case (ii) with $\mu_c=0$. Connecting the minima of $W(z,\cdot)$ at
$\alpha=\AA=0$, at $\alpha=\BB$ and at $\alpha=2\pi$ with minimal energy $0$
yields $W^{**}(z,\alpha)\equiv0$.
Bottom: The case (iii) with $0<\mu_c\le\muccrit$ and $\mu\not=\mu_c$.
For $\alpha\in[0,\AA)$, $W(z,\cdot)$ is strictly convex such that
$W(z,\alpha)=W^{**}(z,\alpha)$. For
$\alpha\in[\AA,\BB]$, $W^{**}(z,\cdot)\equiv\emin-\frac\mu2z^2$ is constant,
connecting the two minima by a straight line. For $\alpha\in[\BB,2\pi]$,
$W^{**}(z,\alpha)$ is a slightly increasing linear function, connecting
$(\BB,W(z,\BB))$ with $(2\pi,W(z,2\pi))$.
In both plots, the double-well is strongly exaggerated to better
illustrate the principle.}
\end{figure}

\vspace*{2mm}
In the most involved case (iii) with $0<\mu_c\le\muccrit$, $\mu\not=\mu_c$,
$W(z,\cdot)$ is convex in $[0,\AA)$, leading to (\ref{Qss3})$_1$.
For $\alpha\in[\AA,\BB]$, $W(z,\cdot)$ forms a double-well potential.
Connecting the minima $(\AA,W(z,\AA))=(\AA,\emin\!-\!\frac{\mu}{2}z^2)$ and
$(\BB,W(z,\BB))=(\BB,\emin\!-\!\frac{\mu}{2}z^2)$ in a straight line leads to
Eqn.~(\ref{Qss3})$_2$.
For $\alpha\in[\BB,2\pi)$, the convexification of $W(z,\cdot)$ is given
by a straight line connecting $(\BB,W(z,\BB))=(\BB,\emin\!-\!\frac{\mu}{2}z^2)$
and $(2\pi,W(z,2\pi))=(2\pi,W(z,0))=(2\pi,\frac{\mu_c}{2}z^2)$, yielding
(\ref{Qss3})$_3$. Fig.~\ref{fig6} illustrates the construction.

\vspace*{2mm}
The construction in the case (i) with $\mu=\mu_c$ or $\mu_c>\muccrit$ is
similar. For $\alpha\in[0,\alpha_2)$,
$W(z,\cdot)$ is convex, leading to Eqn.~(\ref{Qss1})$_1$. For
$\alpha\in[\alpha_2,2\pi)$, the points
$(\alpha_2,W(z,\alpha_2))$
and $(2\pi,W(z,2\pi))$
are connected by a straight line, leading to Eqn.~(\ref{Qss1})$_2$.

For the case (ii) with $\mu_c=0$, connecting the three points
$(\AA,W(z,\AA))=(0,0)$,
$(\BB,W(z,\BB))=\big(\arctan\!\big(\frac{4\gamma}{4-\gamma^2}\big),0\big)$
and $(2\pi,W(z,2\pi))=(2\pi,0)$ where $W(z,\cdot)$ has minimal energy $\emin=0$
yields at once $W^{**}(z,\alpha)\equiv0$. This implies (\ref{Qss2}).
Fig.~\ref{fig6} sketches the construction of $W^{**}(z,\cdot)$ when $\mu_c=0$.
\qed

\newpage
\section{First order Gamma-limit}
\label{secfirst}
For small fixed $\ve>0$,
two minimizers of
$E_\ve(\uhom,\cdot)$ are present, either when $\mu>\mu_c$, $0<\mu_c\le\muccrit$
or when $\mu_c=0$, giving rise to transition layers.
In these cases where there is ambiguity which minimizer is selected, the
zero-order $\Gamma$-limit does not provide the complete physical picture.
For that reason, we may investigate the problem further by studying the
first-order term in the $\Gamma$-expansion of $E$ with scaling $\ve^{-1}$.
It turns out that this rescaled limit is far better suited to reveal the fine
properties of the optimal micro-rotations $\alpha$ and the associated
transition layers than the zero-order $\Gamma$-limit computed in
Proposition~\ref{prop1}, in particular in the case $\mu_c=0$ where
$E_0(\uhom,\cdot)$ is constant over the entire range of $\alpha$ values.

These considerations motivate for $\ve\searrow0$ the introduction of the
rescaled energy
\begin{equation}
\label{E1eps}
F_\ve(u,\alpha):=\frac{E_\ve(u,\alpha)-\emin}{\ve}.
\end{equation}
We will utilize the well-developed machinery of $\Gamma$-expansion outlined in
\cite{BT08,Braides02}. The analytical tools of this section, especially in
Section~\ref{secsurf}, have been developed and sharpened in several articles
on phase transition phenomena, starting with the seminal paper \cite{Modica87}
based on Geometric Measure Theory.

\vspace*{2mm}
We split $\emin$ (cf. Table~\ref{tab1}) in one component for $W$ and one
component due to $\frac{\mu}{2}|\uhom'|^2$ and write for the minima of $W$
\begin{equation}
\label{split}
\wmin:=\emin-\frac{\mu}{2}\gamma^2.
\end{equation}
The splitting (\ref{split}) motivates the introduction of the shifted energies
\begin{align}
\label{V1def}
V_1(u') \;:=\;& \frac{\mu}{2}\Big||u'|^2-\gamma^2\Big|,\\
\label{V2def}
V_2(u',\alpha) \;:=\;& W(u',\alpha)-\wmin\\
\;=\;& \frac{\mu}{2}\Big(\sa u'\!-\!4\shalf\!\Big)^2
+\frac{\mu_c}{2}\Big(\ca u'\!-\!2\sa\Big)^2-\wmin.\nn
\end{align}
The modulus $|\cdot|$ in (\ref{V1def}) ensures the (local) non-negativity of
$V_1$. This prevents a tradeoff between $V_1$ and $V_2$ in such a way that
$V_1(u')$ is locally negative as $|u'|^2<\gamma^2$ in an interval
$I\subset\Omega$ while $V_2(u',\alpha)$ is positive in $I$.
The definition (\ref{V2def}) implies $V_2\ge0$. The zero set of
$V_2(\uhom',\cdot)$ is investigated below in Section~\ref{secgam},
Eqn.~(\ref{zerodef}).

Eqn.~(\ref{E1eps}) and the splitting (\ref{split}) give rise to the definition
\begin{equation}
\label{1order}
F_\ve(u,\alpha):=\left\{ \begin{array}{ll}
\frac{1}{\ve}\int_0^1\ve^2|\alpha'|^2+V_1(u')+V_2(u',\alpha)\dx, &
\mbox{ if }(u,\alpha)\in{\cal X},\\
+\infty, & \mbox{ else.} \end{array}\right.
\end{equation}
\vspace*{2mm}

\subsection{Surface energy}
\label{secsurf}
As can be learned from Table~\ref{tab1} and Figs.~\ref{fig4} and \ref{fig5},
in the regime $\mu_c=0$ and in the regime $\mu>\mu_c$, $0<\mu_c\le\muccrit$,
two optimal micro-rotations $\AA$, $\BB$ occur.
In general, let $\alpha^-$ and $\alpha^+$ be two locally minimizing rotations.
The {\em surface energy} or {\em interfacial energy} due to a transition between
adjacent sets $\{\alpha=\alpha^-\}$ and $\{\alpha=\alpha^+\}$ in $\Omega$ is
\begin{align}
\label{c0def}
c_0(\alpha^-,\alpha^+) \;:=&\; 2\,\Bigg|\!\int_{\alpha^-}^{\alpha^+}
V_2(\uhom',x)^{1/2}\dx\Bigg| \;=\; 2\,\Bigg|\!\int_{\alpha^-}^{\alpha^+}
V_2(\gamma,x)^{1/2}\dx\Bigg|.
\end{align}
Before continuing, we want to comment on this formula. Let $(\ve_n)_{n\in\N}$
be a sequence of positive real numbers with $\ve_n\searrow0$ for $n\to\infty$
and let $(\uen,\aen)_{n\in\N}\subset{\cal X}$ be a minimizing sequence of
$F_{\ve_n}$ as $n\to\infty$.
As will be shown in Lemma~\ref{lem3}, $\aen\to\alpha$ in
measure for a piecewise constant function $\alpha\in L^\infty(\Omega)$. Let us
assume that $\alpha$ jumps at $\overline{x}\in\Omega$, i.e. for a small $t>0$,
$\aen(\overline{x}-t)\to\alpha^-$ and $\aen(\overline{x}+t)\to\alpha^+$
for $n\to\infty$. If we plug in the limit $\alpha$ into $F_{\ve_n}$, we obtain
\begin{equation}
\label{Fs}
F_{\ve_n}(\uen,\alpha)\ge\int_{\overline{x}-t}^{\overline{x}+t}\ve_n
|\alpha'(x)|^2+\frac{1}{\ve_n}V_2(\uen'\!(x),\alpha(x))\dx.
\end{equation}
The following rescaling is a modification of an argument by Modica and Mortula,
\cite{MM77}. After substituting the stretched spatial coordinate
$y:=\frac{x-\overline{x}}{\ve_n}$ and with the functions
\[ \ta(y):=\alpha(\ve_ny+\overline{x}),
\qquad \tu\!(y):=u_{\ve_n}\!\Big(\frac{\ve_ny+\overline{x}}{\ve_n}\Big), \]
we find
\begin{align*}
\ta'(y) &= \frac{\d\ta}{\dy}=\frac{\d\ta}{\d x}\frac{\d x}{\d y}
=\frac{\d\ta}{\d x}\ve_n=\alpha'(x)\ve_n,\\
u_{\ve_n}'\!(x) &= \frac{\d}{\d x}u_{\ve_n}(\ve_ny+\overline{x})
=\frac{\d}{\d y}\tu\!(\ve_ny)\frac{\d y}{\d x}
=\ve_n\tu'\!(\ve_ny)\frac{1}{\ve_n}=\tu'\!(\ve_ny).
\end{align*}
Setting $T_n:=T(\ve_n):=\frac{t}{\ve_n}$, Eqn.~(\ref{Fs}) becomes
\begin{equation}
\label{CS2}
F_{\ve_n}(\uen,\alpha)\ge\int_{-T_n}^{+T_n}|\ta'(y)|^2
+V_2(\tu'\!(\ve_ny),\ta(y))\dy.
\end{equation}
This leads directly to the lower bound
\begin{align}
\hspace*{-3pt}F_{\ve_n}(\uen,\alpha)\,\ge\,
\inf\!\Bigg\{\!\!\int_{-T_n}^{T_n}\!|\ta'(y)|^2
\!+\!V_2(\tu'\!(\ve_ny),\ta(y))\dy\,\Bigg|& \,\ta\!\in\! W^{1,2}(-T_n,T_n);
\nn\\[-20pt]
\label{LB}
& \hspace*{3pt}\ta(-T_n)\!=\!\alpha^-,\,\ta(T_n)\!=\!\alpha^+\!\Bigg\}.
\end{align}
Using the Cauchy-Schwarz inequality $a^2+b^2\ge2\,ab$ yields
\begin{equation}
\label{CS}
\int_{-T_n}^{T_n}\!\!|\ta'(y)|^2+V_2(\tu'\!(\ve_ny),\ta(y))\dy\;\ge\; 2\,
\Bigg|\!\int_{-T_n}^{T_n}\ta'(y)\,V_2(\tu'\!(\ve_ny),\ta(y))^{1/2}\dy\Bigg|.
\end{equation}
In Lemma~\ref{lem3} we will show that $u_{\ve_n}\to\uhom$ in measure for
$n\to\infty$. Consequently, for a sub-sequence, the right hand side of
(\ref{CS}) converges in the limit $n\to\infty$ to
\begin{equation}
\label{cdef}
2\,\Bigg|\!\int_{-\infty}^{\infty}\ta'(y)V_2(\gamma,\ta(y))^{1/2}\dy\Bigg|
=\; 2\,\Bigg|\!\int_{\alpha^-}^{\alpha^+}V_2(\gamma,s)^{1/2}\ds\Bigg| \;=\;
c_0(\alpha^-,\alpha^+).
\end{equation}
Here we substituted $s:=\ta(y)$.

\vspace*{2mm}
This shows that $2\big|\int_{\alpha^-}^{\alpha^+}V_2(\gamma,s)^{1/2}\ds\big|$
in (\ref{cdef}) provides a lower bound on the surface energy which is also
useful for the analysis of the $\Gamma$-limit.
We still have to show that the lower bound is attained, e.g. we need to find
the optimal profile $\ta$.
Following \cite{Sternberg88}, the optimal $\ao$ is available as the solution
to an ordinary differential equation. Let
\begin{align}
\label{ODE1}
\ao'(y) \;&=\; V_2(\gamma,\ao(y))^{1/2}\qquad \mbox{for }y\in\R,\\
\label{ODE2}
\ao(0) \;&=\; \frac{\alpha^-+\alpha^+}{2}.
\end{align}
The existence of a unique solution to (\ref{ODE1}), (\ref{ODE2}) is guaranteed
by the Lipschitz continuity of $\alpha\mapsto V_2(\gamma,\alpha)^{1/2}$ and the
Picard-Lindel{\"o}f theorem. The reason for choosing (\ref{ODE1}), (\ref{ODE2})
is that it is optimal in (\ref{CS}) for $n\to\infty$ since $\ao'(y)>0$ and
\[ |\ao'(y)|^2+V_2(\gamma,\ao(y))=2V_2(\gamma,\ao(y))^{1/2}\ao'(y). \]
In (\ref{ODE1}), it holds $\lim_{y\to\infty}\ao(y)=\alpha^+$,
$\lim_{y\to-\infty}\ao(y)=\alpha^-$. Hence, while for
large $\ve_n$ the profile $\ao$ needs to be slightly modified to fulfil the
boundary conditions $\ta(-T_n)\!=\!\alpha^-$, $\ta(T_n)\!=\!\alpha^+$
to be valid in (\ref{LB}), in the limit $n\to\infty$, $\ao$ is feasible in
the minimization (\ref{LB}). As $\ao$ satisfies (\ref{CS}) with equality in the
limit $n\to\infty$, this demonstrates that $\ao$ is the minimizer in (\ref{LB})
and concludes the derivation of (\ref{c0def}).

{\renewcommand{\arraystretch}{1.15}
\begin{table}[h!bt]
\hspace*{150pt}
\begin{tabular}{|r|l|l|l|}\hline
\multicolumn{1}{|c|}{\mbox{\footnotesize $\gamma$}}&
\multicolumn{1}{|c|}{\mbox{\footnotesize $\BB$}}&
\multicolumn{1}{|c|}{\mbox{\footnotesize $c_0$}}&
\multicolumn{1}{|c|}{\mbox{\footnotesize $\cored$}}\\
\hline\hline
$0.1$ & $0.09917$ & $0.000332$ & $0.000334$\\
$0.2$ & $0.19934$ & $0.00265$  & $0.002666$\\
$0.3$ & $0.29778$ & $0.00888$  & $0.009$\\
$0.4$ & $0.39479$ & $0.020836$ & $0.021334$\\
$0.5$ & $0.48996$ & $0.04017$  & $0.041666$\\
$0.6$ & $0.58291$ & $0.068346$ & $0.072$\\
$0.7$ & $0.67335$ & $0.106602$ & $0.114334$\\
$0.8$ & $0.76101$ & $0.155948$ & $0.170666$\\
$0.9$ & $0.84571$ & $0.217168$ & $0.243$\\
$1.0$ & $0.9273$  & $0.29082$  & $0.333334$\\
\hline
\end{tabular}
\caption{\label{tab2}
The case $\mu_c=0$: Values of $\BB$, $c_0=c_0(\AA,\BB)$ and $\cored$ as a
function of $\gamma$ for $\mu=2$. (Generically, it holds $\AA=0$ for
$\mu_c=0$.)}
\end{table}}

\vspace*{2mm}
For illustration of 
(\ref{c0def}), let us address the case $\mu_c=0$. From (\ref{alphadef}) we have
$\AA=0$, $\BB=\arctan\!\big(\frac{4\gamma}{4-\gamma^2}\!\big)>0$, $\wmin=0$,
leading to
\begin{align*}
c_0 \;:=\; c_0(\AA,\BB) \;=&\; 2\int_{\AA}^{\BB}V_2(\uhom',x)^{1/2}\dx \;=\;
2\Big(\frac{\mu}{2}\Big)^{1/2}\int_{\AA}^{\BB}\Big|\gamma\sin(x)
-4\sin^2\!\big(\frac{x}{2}\big)\Big|\dx\\
=&\; \big(2\mu\big)^{1/2}\int_0^{\BB}\gamma\sin(x)+2\cos(x)-2\dx\\
=&\; 
\big(2\mu\big)^{1/2}\Big[\gamma(1-\cos(\BB))+2\sin(\BB)-2\BB\Big].
\end{align*}
Table~\ref{tab2} lists $c_0$ for $\mu=2$ and different values of $\gamma$.
As a comparison, the table also shows the (numerically computed) values of the
corresponding surface energy $\cored$ for the reduced energy $\Ered$ given
by (\ref{Ered}).

\vspace*{3mm}
\subsection{Computation of the first-order $\Gamma$-limit}
\label{secgam}
We introduce the {\it zero-set} of $\alpha\mapsto V_2(\uhom,\alpha)$ as
\begin{equation}
\label{zerodef}
\Zmu:=\left\{ \begin{array}{ll} \{\AA,\BB\} & \mbox{ if }\mu_c=0\mbox{ or }
\big(\mu>\mu_c\mbox{ and }0<\mu_c\le\muccrit\big),\\
\{\alpha_2\} & \mbox{ if }\mu=\mu_c\mbox{ or }\big(\mu_c>\muccrit\mbox{ and }
\mu\not=\mu_c\big). \end{array}\right.
\end{equation}

Let ${\mathrm PC}(\Omega)\subset L^\infty(\Omega)$ denote the space of
{\it piecewise constant functions} in $\Omega$.
By $\Smu$ we denote the {\it jump set} of $\alpha\in{\mathrm PC}(\Omega)$ with
$\alpha\in\Zmu$ a.e. in $\Omega$.
If $\Zmu$ contains only one element, then $\Smu=\emptyset$.

\begin{Lemma}[Equi-Coercivity of $F_{\ve_n}$]
\label{lem3}
Let $(\ve_n)_{n\in\N}$ be a sequence of positive numbers with $\ve_n\searrow0$
for $n\to\infty$. Let $(u_n,\alpha_n)_{n\in\N}\subset (L^1(\Omega))^2$ with
$\sup_{n\in\N}F_{\ve_n}(u_n,\alpha_n)<\infty$. Then there exists a subsequence
of $(u_n,\alpha_n)$ that converges to $(\uhom,\alpha)$ in $L^1(\Omega)$ with
$\alpha\!\in\!\Zmu$ a.e.
\end{Lemma}

\vspace*{2mm}
\noindent{\bf Proof}. (i) {\it Convergence of $u_n\to\uhom$, $\alpha_n\to\Zmu$
in measure.}

For $n\in\N$ and fixed $\delta>0$, let
\[ I_n^\delta:=\big\{x\in\Omega\;\big|\;{\rm dist}(\alpha_n(x),\Zmu)>\delta
\mbox{ and }|u_n(x)-\uhom(x)|>\delta\big\}. \]
Then it holds for every $n\in\N$
\begin{align*}
& \hspace*{-10pt}
|I_n^\delta|\min_{x\in I_n^\delta}\!\Big(
V_1(u_n'(x))\!+\!V_2(u_n'(x),\alpha_n(x))\Big)\;=\;
\min_{x\in I_n^\delta}\!\Big(V_1(u_n'(x))\!+\!V_2(u_n'(x),
\alpha_n(x))\Big)\int_{I_n^\delta}1\dx\nn\\
\hspace*{25pt} &\le\; \int_{I_n^\delta}V_1(u_n'(x))+V_2(u_n'(x),\alpha_n(x))\dx
\; \le\; \ve_n F_{\ve_n}(u_n,\alpha_n)\le C\ve_n.
\end{align*}
Hence, $|I_n^\delta|\to0$ as $n\to\infty$ for each $\delta>0$.

\vspace*{2mm}
\noindent(ii) {\it $L^1$-convergence of a subsequence.}

The convergence of $(u_n,\alpha_n)\to(\uhom,\alpha)$ in measure implies the
almost everywhere convergence of a subsequence to $(\uhom,\alpha)$ in $\Omega$,
see, e.g., \cite{Folland99}.
In addition, for $\alpha_n\to\alpha$ a.e. in $\Omega$, there exists a further
subsequence $(\alpha_{n_k})_{k\in\N}$ with $\alpha_{n_k}\to\alpha$ in
$L^1(\Omega)$. (For a proof, pick $\alpha_{n_k}$ with
$\|\alpha_{n_{k+1}}-\alpha_{n_k}\|_{L^1(\Omega)}<2^{-k}$ such that
$\sum_{k=1}^\infty \|\alpha_{n_{k+1}}-\alpha_{n_k}\|_{L^1}<\infty$ and use the
monotone convergence theorem.)
\qed

\vspace*{2mm}
As a consequence of Lemma~\ref{lem3}, the fundamental Theorem~\ref{theo1}
applies to $(F_{\ve_n})_{n\in\N}$ and the limit functional $F_0$.

\vspace*{2mm}
We are now ready to prove the following main result.

\begin{Proposition}
\label{prop2}
Let $F_\ve$ be defined by (\ref{1order}) and let $c_0(\alpha^-,\alpha^+)$
be given by (\ref{c0def}).\\
Then the $\Gamma-lim_{\ve\searrow0}F_\ve$ with respect to convergence in
$L^1(\Omega)$ exists and is given by
\begin{equation}
\label{limi}
F_0(u,\alpha)=\left\{\!\begin{array}{ll}\sum\limits_{\Smu}
c_0(\alpha^-,\alpha^+), & \mbox{if }u=\uhom,\, u(0)\!=\!0,\,u(1)\!=\!\gamma,
\\[-2mm]
& \quad \alpha\!\in\!{\mathrm PC}(\Omega),\,\alpha\in\Zmu\mbox{ a.e.}\\
+\infty, & \mbox{else.}\end{array}\right.
\end{equation}
Especially, $F_0(u,\alpha)=0$ for $\Smu=\emptyset$.
\end{Proposition}

\newpage
\noindent{\bf Proof.} (i) {\it Liminf inequality}.

Let $(\ve_n)_{n\in\N}$ be a sequence of positive real numbers with
$\ve_n\searrow0$ as $n\to\infty$.
Let $(\uen,\aen)\to(u,\alpha)$ in $L^1(\Omega)$ for $n\to\infty$
and $\sup_{n\in\N}F_{\ve_n}(\uen,\aen)<\infty$. We have to show that
\begin{equation}
\label{ls}
F_0(u,\alpha)\le\liminf_{n\to\infty}F_{\ve_n}(\uen,\aen).
\end{equation}
If $\Smu=\emptyset$, then $F_0(u,\alpha)=0$ and
(\ref{ls}) is evidently true.

\vspace*{2mm}
By the equi-coercivity of $F_{\ve_n}$ proven in Lemma~\ref{lem3}, we already
know that $u=\uhom$ a.e. in $\Omega$ and $\alpha\in{\mathrm PC}(\Omega)$ with
$\alpha\in\Zmu$ a.e. in $\Omega$. For $\Smu\not=\emptyset$, we first consider
the case that
$\alpha$ possesses exactly one jump point $\overline{x}\in\Omega$, i.e. there
is a $t>0$ with $\aen(\overline{x}-t)\to\alpha^-$,
$\aen(\overline{x}+t)\to\alpha^+$ as $n\to\infty$ with
real values $\alpha^-\not=\alpha^+$. Using (\ref{CS2}) and (\ref{CS}), for
$T_n:=\frac{t}{\ve_n}$,
\begin{align}
F_{\ve_n}(\uen,\aen) \;\ge\; & \int_{-T_n}^{T_n}\ve_n|\ta_{\ve_n}'(y)|^2
+\frac{1}{\ve_n}V_2(\tu'\!(\ve_ny),\ta_{\ve_n}(y))\dy\nn\\
\label{linf1}
\;\ge\; & 2\,\Bigg|\!\int_{-T_n}^{T_n}\ta_{\ve_n}'(y)V_2(\tu'\!(\ve_ny),
\ta_{\ve_n}(y))^{1/2}\dy\Bigg|.
\end{align}
The right hand side of (\ref{linf1}) converges for $n\to\infty$ to
\[ 2\,\Bigg|\!\int_{-\infty}^{\infty}\ta'(y)V_2(\gamma,\ta(y))^{1/2}\dy\Bigg|
=2\,\Bigg|\!\int_{\alpha^-}^{\alpha^+}V_2(\gamma,s)^{1/2}\ds\Bigg|=
c_0(\alpha^-,\alpha^+)=F_0(u,\alpha) \]
and we obtain (\ref{ls}) as desired.

If $\Smu$ contains several jump points, let
$\overline{\Omega}=\cup_{i=1}^N[x_i,x_{i+1}]$ with $x_1\!=\!0$ and
$x_{N+1}\!=\!1$ such that each sub-interval $(x_i,x_{i+1})$ contains at most
one element of $\Smu$ and such that $\aen(x_i)\to\alpha(x_i)$ for
$2\le i\le N$. Like in (\ref{linf1}) we end up with
\begin{align}
F_{\ve_n}(\uen,\aen)
\label{linf2}
&\ge \;\sum_{i=1}^N2\,\Bigg|\!\int_{-T_n^i}^{T_n^i}\ta_{\ve_n}'\!(y)
V_2(\tu'\!(\ve_ny,\ta_{\ve_n}\!(y))^{1/2}\dy\Bigg|.
\end{align}
As in the case with only one jump point, for $n\to\infty$ the right hand side
of (\ref{linf2}) converges to
$\sum_{\Smu}c_0(\alpha^-,\alpha^+)\!=\!F_0(u,\alpha)$, proving (\ref{ls}).

\vspace*{2mm}
\noindent(ii) {\it Limsup inequality}.

We need to find a recovery sequence $(u_\ve,\alpha_\ve)\subset (L^1(\Omega))^2$
with
\begin{equation}
\label{rs}
F_0(u,\alpha)\ge\limsup_{\ve\searrow0}F_\ve(u_\ve,\alpha_\ve).
\end{equation}
If $u\not=\uhom$, then $F_0(u,\alpha)=+\infty$ and (\ref{rs}) is obviously true.
Due to the additivity of the integral similar to the reasoning in (i), it is
enough to restrict to the case
\[ \alpha(x)=\left\{ \begin{array}{ll}\alpha^-, & \mbox{ if }x<\overline{x},\\
\alpha^+, & \mbox{ if }x\ge\overline{x} \end{array}\right. \]
where $\overline{x}\in\Omega$ is fixed. When $\Smu=\emptyset$, we have
$\alpha^-=\alpha^+$ and $c_0(\alpha^-,\alpha^+)=0$.

\vspace*{2mm}
Now let $\Smu=\{\overline{x}\}$ and $\alpha^-\not=\alpha^+$.
Fix $\delta>0$. Due to the considerations in Section~\ref{secsurf},
there exists a $T>0$ and a function $\beta\in W^{1,2}(-T,T)$ such that
$\beta(-T)=\alpha^-$, $\beta(T)=\alpha^+$ and
\[ \int_{-T}^T|\beta'(y)|^2+V_2(\uhom',\beta(y))\dy\;\le\;c_0(\alpha^-,\alpha^+)
+\delta. \]
We can construct a recovery sequence by choosing
\begin{equation}
\label{star2}
\alpha_\ve(x):=\left\{\begin{array}{ll} \alpha^-, & \mbox{ if }
x<\overline{x}-\ve T,\\
\beta(\ve x+\overline{x}), & \mbox{ if }\overline{x}-\ve T\le x\le
\overline{x}+\ve T,\\
\alpha^+, & \mbox{ if }x>\overline{x}+\ve T. \end{array}\right.
\end{equation}
This yields
\begin{align*}
\limsup_{\ve\searrow0}F_\ve(u_\ve,\alpha_\ve) \; &= \; \int_{-T}^T|\beta'(y)|^2
+V_2(\gamma,\beta(y))\dy\\
&\le\; c_0(\alpha^-,\alpha^+)+\delta\;=\;F_0(u,\alpha)+\delta.
\end{align*}
As $\delta>0$ is arbitrary, in the limit $\delta\searrow0$ we end up
with (\ref{rs}).

If $\Smu$ contains several jump points, we observe that by (\ref{star2}),
$\alpha$ is only modified on a small set near $\Smu$. Hence, the above
construction can also be carried out for an arbitrary
$\alpha\in{\mathrm PC}(\Omega)$ with $\alpha\in\Zmu$ a.e. \qed

\section{Conclusion}
\label{seccon}
In this paper, the zero and first order Gamma limit of $E(u,\alpha)$ have been
computed and the minimizers have been identified. In particular, the results
reveal the fine properties of the optimal micro-rotations $\alpha$ forming
transition layers in $\Omega$.

The relaxed functionals $E_0$ and $F_0$ may also be of interest
for numerical simulations. Using $E_0$ instead of the original Cosserat
functional $E$ given by (\ref{Edef}) for simulations with a small but finite
$L_c>0$ corresponds to a convexification or homogenization of the problem
and may help apart from a very significant speed up to avoid some of the
numerical problems encountered in \cite{Ble15,Ble17,Ble20}.

\newpage
\section*{Appendix -- Direct minimization of $E$ for $L_c=0$}
\label{secapp}
It is instructive to compare the zero-order Gamma-limit of $E_\ve$, i.e. the
results of Proposition~\ref{prop1} and Lemma~\ref{lem2}, with the following
{\it direct minimization}. Let $\widetilde{E}_0$ be defined by (\ref{E0tdef})
which coincides with $E$ given by Eqn.~(\ref{Edef}) with $L_c=0$.
For chosen deformation $u\in{\cal X}_u$, we denote by
$\aopt=\aopt(u')$ a corresponding optimal micro-rotation, i.e. a rotation
$\alpha\in{\cal X}_a$ that minimizes $\alpha\mapsto E(u,\alpha)$ for fixed $u$.
Plugging in $\aopt$ into $\widetilde{E}_0$, we end up with the functional
\begin{align}
E_\mathrm{opt}(u)\; :=\; & \widetilde{E}_0(u,\aopt(u'))\nn\\
=\; & \frac\mu2 \int_0^1|u'|^2+\big(\sin(\aopt)u'-4\sin^2\!\big(
\frac{\aopt}{2}\!\big)\big)^2\dx\nn\\
\label{Eoptdef}
& +\frac{\mu_c}{2}\int_0^1\big(\cos(\aopt)u'-2\sin(\aopt)\big)^2\dx.
\end{align}

The following proposition computes $E_\mathrm{opt}$ explicitly for the
different regimes.

\begin{Proposition}[Direct minimization of $E$ for $L_c=0$]
\label{prop3}
Let $E$ be given by (\ref{Eoptdef}) and assume
$u\in{\cal X}_u\cap W^{2,2}(\Omega;\,\R)$. Then it holds

\noindent(i) If $\mu=\mu_c$ or $\mu_c>\muccrit:$
\begin{equation}
\label{Eopt1}
E_\mathrm{opt}(u)=\mu\int_0^1|u'|^2+4-2\big(|u'|^2+4)^{1/2}\dx.
\end{equation}
(ii) If $\mu_c=0:$
\begin{equation}
\label{Eopt2}
E_\mathrm{opt}(u)=\frac\mu2\int_0^1|u'|^2\dx.
\end{equation}
(iii) If $0<\mu_c\le\muccrit$, $\mu\not=\mu_c:$
\begin{equation}
\label{Eopt3}
E_\mathrm{opt}(u)=\frac{\mu+\mu_c}{2}\int_0^1|u'|^2\dx
-\frac{2\mu_c^2}{\mu\!-\!\mu_c}.
\end{equation}
The functional $E_\mathrm{opt}$ defined by any of (\ref{Eopt1}), (\ref{Eopt2}),
or (\ref{Eopt3}) is convex in $u'$.
\end{Proposition}

\vspace*{2mm}
\noindent{\bf Proof.}
Considering the two Euler-Lagrange equations of $\widetilde{E}_0$
defined in (\ref{E0tdef}), the first constitutes a balance of forces while
the second states an algebraic relationship between $\alpha$ and
$u'$ since $L_c=0$. Resolving this second equation leads to the
minimizing optimal rotations $\aopt=\aopt(u')$ summarized in Table~\ref{tab1}.

With the regularity assumption
$u\in W^{2,2}(\Omega;\,\R)$, it holds $u''=0$ in $\Omega$, see
\cite[Lemma~3]{BN22}. From the Sobolev embedding
$W^{2,2}(\Omega)\hookrightarrow C^1(\overline{\Omega})$ we infer
$u'\equiv\mathrm{const}$ in $\Omega$.

\vspace*{2mm}
\noindent(i) Let $\mu=\mu_c$ or $\mu_c>\muccrit$.
The unique optimal rotation in this case is
$\aopt(u')=\alpha_2=\arctan\!\big(\frac{u'}{2}\big)$. For $t\in\R$ we remark
the identities
\begin{align}
\label{H1}
\cos(\arctan(t)) \;=&\; \frac{1}{(t^2+1)^{1/2}},\\
\label{H2}
\sin(\arctan(t)) \;=&\; \frac{t}{(t^2+1)^{1/2}}.
\end{align}
With (\ref{H1}), (\ref{H2}), direct inspection yields
$\sin(\alpha_2)=\frac{u'}{(|u'|^2+4)^{1/2}}$,
$\cos(\alpha_2)=\frac{2}{(|u'|^2+4)^{1/2}}$ such that
\begin{align*}
\cos(\alpha_2)u'-2\sin(\alpha_2) &= 0,\\ 
\big(\sin(\alpha_2)u'+2\cos(\alpha_2)-2\big)^2 &= \Big(
\frac{|u'|^2+4}{(|u'|^2+4)^{1/2}}-2\Big)^2 \;=\; |u'|^2+8-4\big(|u'|^2+4)^{1/2}.
\end{align*}
Using the identity $-4\sin^2\!\big(\frac{\aopt}{2}\!\big)=2\cos(\aopt)-2$,
this shows for case (i)
\[ E_\mathrm{opt}(u)=\frac{\mu}{2}\int_0^1|u'|^2+\big(\sin(\alpha_2)u'
+2\cos(\alpha_2)-2\big)^2\dx=\mu\int_0^1|u'|^2+4-2\big(|u'|^2+4\big)^{1/2}\dx \]
which is (\ref{Eopt1}).

\vspace*{2mm}
\noindent(ii) Let $\mu_c=0$. The optimal rotations in this case are
\[ \aopt=\AA=0,\qquad \aopt(u')=\BB=\eta^{-1}(u'), \]
cf. Table~\ref{tab1}, where $\eta$ is defined by (\ref{a3def}).
In both cases, $\big(\sin(\aopt)u'-4\sin^2(\aopt/2)\big)^2=0$, and
(\ref{Eoptdef}) at once simplifies to (\ref{Eopt2}).

\vspace*{2mm}
\noindent(iii) Let $0<\mu_c\le\muccrit$, $\mu\not=\mu_c$.
With $f:=\big((|u'|^2+4)(\mu\!-\!\mu_c)^2-4\mu^2\big)^{1/2}$, cf.
Eqn.~(\ref{fdef}), the two optimal rotations in this case are
\[ \aopt(u')=\BB=\arctan\!\Big(\frac{\mu u'+f}{2\mu-\frac{u'}{2}f}\Big),\qquad
\aopt(u')=\AA=\arctan\!\Big(\frac{\mu u'-f}{2\mu+\frac{u'}{2}f}\Big). \]
By direct inspection, we find
\begin{align*}
\sin(\AA)u'+2\cos(\AA)=\frac{2\mu}{\mu-\mu_c},\qquad
\cos(\AA)u'-2\sin(\AA)=\frac{f}{\mu-\mu_c}, \\
\sin(\BB)u'+2\cos(\BB)=\frac{2\mu}{\mu-\mu_c},\qquad
\cos(\BB)u'-2\sin(\BB)=\frac{-f}{\mu-\mu_c}.
\end{align*}
Plugging these identities into (\ref{Eoptdef}),
we obtain for both choices of $\aopt$
\begin{align*}
E_\mathrm{opt}(u) =&
 \frac{\mu}{2}\int_0^1|u'|^2+\Big(\frac{2\mu}{\mu\!-\!\mu_c}-2\Big)^2\dx
+\frac{\mu_c}{2}\int_0^1\frac{f^2}{(\mu\!-\!\mu_c)^2}\dx\\
=& \frac{\mu}{2}\int_0^1|u'|^2\dx
+\frac{2\mu\mu_c^2}{(\mu\!-\!\mu_c)^2}+\frac{\mu_c}{2}\int_0^1
\frac{(|u'|^2+4)(\mu\!-\!\mu_c)^2-4\mu^2}{(\mu\!-\!\mu_c)^2}\dx\\
=& \frac{\mu+\mu_c}{2}\int_0^1|u'|^2\dx+\frac{2\mu\mu_c^2}{(\mu\!-\!\mu_c)^2}
+\frac{\mu_c}{2}\;\frac{4\mu_c^2\!-\!8\mu\mu_c}{(\mu\!-\!\mu_c)^2}.
\end{align*}
This simplifies to (\ref{Eopt3}).

The convexity of $E_\mathrm{opt}$ given by Eqn.~(\ref{Eopt2}) and
(\ref{Eopt3}) is evident. But also Eqn.~(\ref{Eopt1}) defines a convex
functional in $z=u'$, even though it may first not appear so.
Indeed, introducing $g:\R\to\R$,
\[ g(z):=z^2+4-2\big(z^2+4\big)^{1/2}, \]
a direct computation yields $g'(z)=2z-2z\big(z^2+4)^{-1/2}$ and
\[ g''(z)=2+\frac{2z^2}{(z^2+4)^{3/2}}-\frac{2}{(z^2+4)^{1/2}}\;=\;
\frac{2(z^2+4)^{3/2}-8}{(z^2+4)^{3/2}}\;>\;0. \]
This is the convexity of $g$ and hence of $E_\mathrm{opt}$ in $u'$ as defined
by Eqn.~(\ref{Eopt1}). \qed

\begin{Remark}
\label{rem2}
For the case (ii) with $\mu_c=0$, $E_\mathrm{opt}(u)$ coincides with 
$E_0:=\Gamma\!-\!\lim_{\ve\searrow0}E_\ve$ computed in Prop.~\ref{prop2} and
Lemma~\ref{lem2}.
For the other cases (i) and (iii) of Proposition~\ref{prop3}, $E_\mathrm{opt}$
differs from the Gamma-limit $E_0$.
This underlines the critical role of the Cosserat couple modulus $\mu_c$
in the modelling.
\end{Remark}

\begin{Remark}
\label{rem3}
A direct minimization analogous to (\ref{Eoptdef}) is also possible in three
space dimensions for $\widetilde{E}_{3D}$, cf. Eqn.~(\ref{E3ddef}).
In \cite{NFB19,FN17a,FN17b}, the optimal rotations are computed.
However, for $\mu_c\ge\mu$ it is known, \cite{Neff14}, that 
the resulting functional
\[ \int_{\widehat{\Omega}}\mu\,\mathrm{dist}^2\big(F,\mathrm{SO}(n)\big)
+\frac{\lambda}{4}\Big[(\det\overline{U}\!-\!1)^2+\Big(
\frac{1}{\det\overline{U}}\!-\!1\Big)^2\Big]\dx \]
is not rank-one convex due to the $\mathrm{dist}$-function.
The computation of the quasi-convex hull w.r.t. deformations in
$\mathrm{GL}^+(2)$ in this case can be found in \cite{GM23}.
\end{Remark}


\vspace*{2mm}
\begin{thebibliography}{00}
\bibitem{appell1893}
Appell, P. (1893).
Trait\'{e} de m\'{e}canique rationnelle: Statique. Dynamique du point
(Vol. 1). Gauthier-Villars.
\bibitem{BN22}
Blesgen, T., Neff, P. (2022). Simple shear in nonlinear Cosserat micropolar
elasticity: Existence of minimizers, numerical simulations and occurrence of
microstructure, Mathematics and Mechanics of Solids, 
DOI 10.1177/10812865221122191
\bibitem{Ble13} Blesgen, T. (2013).
Deformation patterning in Cosserat plasticity,
Modelling and Simulation in Materials Science and Engineering, 21(3),
35001--35012.
\bibitem{Ble14} Blesgen,\hspace*{-3pt} T. (2014). 
Deformation patterning in three-dimensional large-strain Cosserat plasticity,
Mechanics Research Communications 62, 37--43
\bibitem{Ble15} Blesgen, T. (2015). On rotation deformation zones for
finite-strain Cosserat plasticity, Acta Mechanica 226, 2421--2434.
\bibitem{Ble17} Blesgen, T. (2017). A variational model for dynamic
recrystallization based on Cosserat plasticity, Composites B 115, 236--243
\bibitem{Ble20} 
Blesgen, T., Amendola, A. (2020).
Mathematical analysis of a solution method for finite-strain holonomic
plasticity of Cosserat materials,
Meccanica 55, 621--636
\bibitem{Boyd09} Boyd, S. (2009) Convex optimization,
Cambridge University Press.
\bibitem{Braides02} Braides, A. (2002). Gamma-convergence for Beginners,
Oxford Lecture Series in Mathematics. 22nd Edition.
\bibitem{BT08}
Braides, A. and Truskinovsky, L. (2008).
Asymptotic expansions by $\Gamma$-convergence,
Continuum Mechanics and Thermodynamics 20, 21--62.
\bibitem{Capriz89}
Capriz, G. (1989). Continua with {M}icrostructure. Springer.
\bibitem{Cos1909}
Cosserat, E., Cosserat, F. (1909).  
Th{\'e}orie des corps d{\'e}formables. Appell, Paul. Gauthier-Villars, Paris.
\bibitem{Cos1991}
Cosserat, E., Cosserat, F. (1991).
Note sur la th{\'e}orie de l'action euclidienne.
Appendix in \cite{appell1893}, 557--629.
\bibitem{Dmit09}
Dmitrieva, O., Dondl, P.W., M{\"u}ller, S., Raabe, D. (2009).
Lamination microstructure in shear deformed copper single crystals.
Acta Materialia, 57(12), 3439--3449.
\bibitem{CH86} Elliott, C., Songmu, Z. (1986).
On the Cahn-Hilliard equation, Archive for Rational Mechanics and Analysis
96(4), 339--357.
\bibitem{FN17a} Fischle, A., Neff, P., (2017). The geometrically
nonlinear Cosserat micropolar shear-stretch energy. Part~I: A general parameter
reduction formula and energy-minimizing microrotations in 2D,
Zeitschrift f{\"u}r Angewandte Mathematik und Mechanik 97(7), 828--842.
\bibitem{FN17b} Fischle, A., Neff, P., (2017). The geometrically
nonlinear Cosserat micropolar shear-stretch energy. Part~II: Non-classical
energy-minimizing microrotations in 3D and their computational validation,
Zeitschrift f{\"u}r Angewandte Mathematik und Mechanik 97(7), 843--871.
\bibitem{Folland99} Folland, G.B. (1999).
Real Analysis, Prentice Hall.
\bibitem{GM23} Ghiba, I.-D., Martin, R.J., K{\"o}hler, M., Balzani, D., Neff,
P. (2023). Quasiconvex relaxation of a planar Biot-type energy on
$\mathrm{GL}^+(2)$ versus $\R^{2\times2}$. Analytical and numerical approaches.
In preparation.
\bibitem{Lecca13}
Lecca, P. (2013). Stochastic chemical kinetics.
Biophys. Rev. Vol. 5(4), 323--345.
\bibitem{DalMaso93} 
Maso, G.D. (1993). The direct method in the Calculus of Variations.
In: An Introduction to $\Gamma$-convergence. Progress in Nonlinear Differential
Equations and Their Applications, Vol. 8, Birkh{\"a}user Boston.
\bibitem{Levine05}
Levine, R.D. (2005). Molecular Reaction Dynamics.
Cambridge University Press.
\bibitem{Modica87}
Modica, L. (1987).
The gradient theory of phase transitions and the minimal interface criterion.
Archive for Rational Mechanics and Analysis, 98, 123--142.
\bibitem{MM77}
Modica, L., Mortula, S. (1977). Un esempio di $\Gamma$-convergenza.
Boll. Un. Mat It. B 14, 285--299.
\bibitem{SM98}
M{\"u}ller, S. (1998).
Variational models for microstructure and phase transitions,
Lecture Notes no. 2, Max-Planck-Institute for Mathematics.\\
{https://www.mis.mpg.de/publications/other-series/ln/lecturenote-0298.html}
\bibitem{Neff06}
Neff, P. (2006).
A finite-strain elastic-plastic Cosserat theory for polycrystals with
grain rotations. International Journal of Engineering Science 44(8-9), 574--594.
\bibitem{Neff15}
Neff, P., B{\^\i}rsan, M., Osterbrink, F. (2015).
Existence theorem for geometrically nonlinear Cosserat micropolar model under
uniform convexity requirements. Journal of Elasticity, 121(1), 119--141.
\bibitem{NFB19}
Neff, P., Fischle, A., and Borisov, L. (2019).
Explicit global minimization of the symmetrized Euclidean distance by a
characterization of real matrices with symmetric square.
SIAM Journal on Applied Algebra and Geometry. 3(1). 31--43.
\bibitem{Neff14}
Neff, P., Lankeit, J., Madeo, A. (2014).
On Grioli's minimum property and its relation to Cauchy's polar decomposition.
International Journal of Engineering Science. 80, 209--217.
\bibitem{Neff08}
Neff, P., M{\"u}nch, I. (2008).
Curl bounds Grad on $\SO(3)$.
ESAIM: Control, Optimisation and Calculus of Variations 14(1), 148--159.
\bibitem{Neff09}
Neff, P., M{\"u}nch, I. (2009).
Simple shear in nonlinear Cosserat elasticity: bifurcation and induced
microstructure. 
Continuum Mechanics and Thermodynamics, 21(3), 195--221.
\bibitem{Sternberg88}
Sternberg, P. (1988). The effect of a singular perturbation on nonconvex
variational problems, Archive for Rational Mechanics and Analysis
1988, 209--260.
\end{thebibliography}
\end{document}